\documentclass[11pt,a4paper]{article}
\usepackage[latin1]{inputenc}
\usepackage{a4wide,amssymb,amsmath}
\usepackage{graphicx}
\usepackage{times}
\newtheorem{theorem}{Theorem}
\newtheorem{lemma}{Lemma}
\newtheorem{definition}{Definition}

\def\endlines#1{\medskip\vskip\parskip\nointerlineskip
  \centerline{\vbox{\hrule width 1in}}
  \medskip\nointerlineskip\vskip\parskip
  \centerline{\vbox{\halign{\hfil##\hfil\cr#1\crcr}}}}

\begin{document}
\title{Geometric continuity and compatibility conditions\\for 4-patch surfaces}
\author{Bo I Johansson\footnote{Department of Mathematical Sciences, Chalmers University of Technology} \footnote{Department of Mathematical Sciences, University of Gothenburg}}
\date{}
%\author{Bo I Johansson}
%\date{Department of Mathematical Sciences\\Chalmers University of Technology and G\"oteborg University}

\maketitle
\medskip\smallskip
%%%%%\markboth{\today}{\today}

\centerline{\vbox{\hrule width 2in}}
\begin{abstract}
  When considering regularity of surfaces, it is its geometry that is
  of interest. Thus, the concept of geometric regularity or geometric
  continuity of a specific order is a relevant concept. In this
  paper we discuss necessary and sufficient conditions for a 4-patch
  surface to be geometrically continuous of order one and two or, in
  other words, being tangent plane continuous and curvature continuous
  respectively.  The focus is on the regularity at the point where the
  four patches meet and the compatibility conditions that must
  appear in this case. In this article the compatibility conditions are proved to be
  independent of the patch parametrization, i.e., the compatibility conditions are universal. In the end of the paper these
  results are applied to a specific parametrization such as Bezier representation in order to generalize a 4-patch surface result by Sarraga.
\end{abstract}

\section{Introduction}
In many applications in Computer-Aided Geometric Design (CAGD) and
Computer Graphics a surface is composed of several patches, where a
patch usually is represented by a Bezier polynomial, B-spline or
NURBS. In particular, each patch is as regular as is needed. Thus, when
considering smoothness of a surface such as tangent plane continuity
or curvature continuity, the lack of regularity may only occur
somewhere at a common boundary curve between two or more patches.

Regularity for a surface constituting of two adjacent patches sharing
the same common boundary curve, see Figure~1, is a well studied
problem. General results for such a 2-patch surface in the case of
$G^{1}$ as well as $G^{2}$ continuity, i.e., tangent plane continuity
and curvature continuity respectively, were given by Juergen Kahmann in
a paper from 1983, see \cite{kahmann}. In the same paper he applies
these results to the case of Bezier patches. Other authors such as Degen
\cite{degen}, Liu and Hoschek \cite{liuHoschek}, Liu \cite{liu},
DeRose \cite{derose}, have also treated tangent plane continuity in
the 2-patch case. In the case of curvature continuity of 2-patch
surface we refer to articles by Kiciak \cite{kiciak}, Ye, Liang and
Nowacki \cite{yeLiangNowacki}.

A more complicated situation is regularity of a surface
consisting of four patches where every pair of adjacent patches meet
at a common boundary curve and all the patches intersect at a common
vertex, see Figure~\ref{fig2}. Among the many authors that have
treated regularity problem in this 4-patch surface case are B\'ezier
\cite{bezier}, Sarraga \cite{sarraga}, \cite{sarraga1} and
\cite{sarraga2}, Ye and Nowacki \cite{yeNowacki}. Further references
and an overview can be found in the book by Hoschek and Lasser
\cite{hoschekLasser}. In refered articles so far there have only been considered certain cases of parametric patches, not a general parametrization as is done in this paper. A general approach has also been done by Peters, see the articles by Peters \cite{peters}--\cite{peters3} and by Ye \cite{ye}.

%In the paper \cite{peters1} Peters has treated the same problem.

A general approach to study regularity for a 4-patch surface has been
to restrict the patch parametri\-zation to a certain explicit polynomial
or rational basis. In this paper we consider 4-patch surfaces where
the patches are given by any function of the form $(u,v)\mapsto
\mathbf{r}(u,v)\in \mathbb R^{3},\,\, u,v\in[0,1]$. We present
necessary and sufficient compatibility conditions in order to have
tangent plane continuity and curvature continuity respectively for
´such surfaces. The results we achieve in this paper are independent of
the patch parametrization. Thus, these results are of a general nature
and can be applied to any parametrization.

In the last section of this paper we apply our results by considering the Sarraga case of filling a hole of a $G^1$-surface in such a way that the extended surface preserves the $G^1$-regularity, see Sarraga \cite{sarraga}. Compared to Sarragas result we reduce the bi-degree to (5,5) of the created interior patch. More generally, we give necessary and sufficient conditions in order to create such an interior patch.

\section{Geometric continuity of order 1}
When discussing regularity of a surface our focus is on the geometry
of the surface
and not on its actual parametrization. Thus the notation of geometric
continuity
is the concept used in this context. The lowest order of regularity is
$G^{0}$, which means that the surface is connected. Another way to put
it is to require that its representation is
continuous. The next level of regularity is tangent plane
continuity, denoted by $G^{1}$, which is defined here.

\begin{definition}\label{defG1}
  A continuous surface is said to be tangent plane continuous, denoted
  by $G^{1}$, if every point on the surface has a unique tangent
  plane, which varies continuously on the surface. Such a surface is
  also said to be geometrically continuous of order one.
\end{definition}

Consider Definition~\ref{defG1} in the case of a 2-patch surface
$S$. % consists of two patches with a common boundary as in Figure~1.
Here we use the notation $\mathbf{r}\in C^{1}_{\#}$ for a patch
described by a continuous differential function $(u,v)\mapsto
\mathbf{r}(u,v)\in \mathbb{R}^{3}$ with $0\le u,v\le1$ satisfying
$\mathbf{r}_{u}\times \mathbf{r}_{v}\ne0$ for $0\le u,v\le1$. Let the
two patches be described as $(u,v)\mapsto \mathbf{r}^{(1)}(u,v)$ with
$0\le u,v\le1$, and $(s,t)\mapsto \mathbf{r}^{(2)}(s,t)$ with $0\le
s,t\le 1$.  Suppose further that each patch is regular enough, i.e.,
$\mathbf{r}^{(1)},\mathbf{r}^{(2)}\in C^{1}_{\#}$.  In order for the
surface $S$ to be tangent plane continuous the only points that do not
automatically fulfill the $G^{1}$-condition are those along the common
boundary curve of the two patches, see Figure~1. On this boundary
curve we must particularly have $v\mapsto
\mathbf{r}^{(1)}(1,v)=\mathbf{r}^{(2)}(0,t(v))$ for $0\le v\le1$,
where $v\mapsto t(v)$ is a regular reparametrization of the interval
$[0,1]$.  \vspace{5mm}

\begin{figure}[hbtp]\label{figge1}
  \begin{center}
    \leavevmode
     \includegraphics[width=6cm]{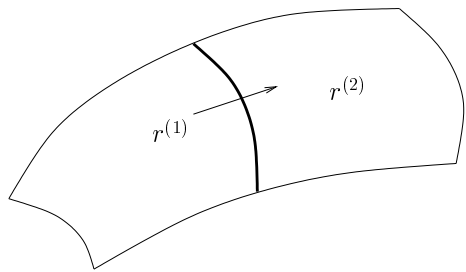}
   \caption{Two patches connected by a common boundary curve}
  \end{center}
\end{figure}

Thus, at a particular boundary point
with parameter value $v$ the tangent plane must fulfill
$$
\mbox{\em span}\{\mathbf{r}^{(1)}_{u}(1,v),\mathbf{r}^{(1)}_{v}(1,v)\} =
\mbox{\em span}\{\mathbf{r}^{(2)}_{s}(0,t(v)),\mathbf{r}^{(2)}_{t}(0,t(v))\}.
$$
Let us formulate this in an alternative way, where we use the notation $\mathbf{r}^{(2)}(u,v)=\mathbf{r}^{(2)}(u,t(v))$ for $0\le u,v\le1$. The statement can be found in e.\ g.\ \cite{yeLiangNowacki} and is summarized in the next Lemma.

 \begin{lemma}\label{lemma1}
   A necessary and sufficient condition for two adjacent $C^{1}_{\#}$-patches\ $\mathbf{r}^{(1)}$ and $\mathbf{r}^{(2)}$ joining $G^1$-continuously along its common boundary curve $v\mapsto
 \mathbf{r}^{(1)}(1,v)=\mathbf{r}^{(2)}(0,v)$ is that there exist continuous functions $\lambda_{1,2}$ and $\kappa_{1,2}$ such that
 \begin{equation}\label{firstEq}
   \mathbf{r}^{(2)}_u(0,v) =
 \lambda_{1,2}(v)\,\mathbf{r}^{(1)}_u(1,v)+\kappa_{1,2}(v)\,\mathbf{r}^{(1)}_v(1,v),\,\,\,\,0\le
 v\le1.
 \end{equation}
 \end{lemma}

\begin{figure}[hbtp]\label{fig2}
  \begin{center}
    \leavevmode
    \includegraphics[width=5cm]{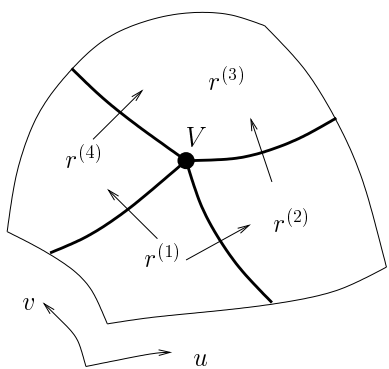}
    \caption{Four patches connected in a common vertex $V$}
  \end{center}
\end{figure}

In this paper we will consider the problem of
a surface constituting of four patches
where every two adjacent patches have a common boundary curve, see
Figure~\ref{fig2}. Moreover, the four patches intersect at a common vertex
$V$. In the case of a 4-patch surface being $G^1$ we will prove that there must exist
compatibility conditions at the intersection point $V$.
% that must be satisfied in order for every two
% patches with a common boundary to fulfill equation (\ref{firstEq}).
% Thus, we want to find the correct conditions for a 4-patch surface to
% be tangentially continuous at the intersection point $V$.
The compatibility conditions can be rephrased in such a way that we
formulate necessary and sufficient conditions on the functions
$\lambda_{ij}$ and $\kappa_{ij}$ at the intersection point $V$.

First, by using
the relation (\ref{firstEq}) we get the
next four relations between the patches $(1)$--$(2)$, $(2)$--$(3)$,
$(4)$--$(3)$ and $(1)$--$(4)$. We use the same parameters $u$ and $v$ for
all the patches, where $0\le u,v\le1$. Thus
\begin{equation}
\begin{split}
\mathbf{r}^{(2)}_u(0,v) & =
\lambda_{1,2}(v)\mathbf{r}^{(1)}_u(1,v)+\kappa_{1,2}(v)\mathbf{r}^{(1)}_v(1,v)\\
\mathbf{r}^{(3)}_u (0,v) & = \lambda_{4,3}(v)\mathbf{r}^{(4)}_u
  (1,v)+\kappa_{4,3}(v)\mathbf{r}^{(4)}_v(1,v)
\end{split}\label{1.1}
\end{equation}
and
\begin{equation}
\begin{split}
\mathbf{r}^{(4)}_v(u,0) & =
\lambda_{1,4}(u)\mathbf{r}^{(1)}_v(u,1)+\kappa_{1,4}(u)\mathbf{r}^{(1)}_u(u,1)\\
\mathbf{r}^{(3)}_v(u,0) & = \lambda_{2,3}(u)\mathbf{r}^{(2)}_v(u,1)+ \kappa_{2,3}(u)\mathbf{r}^{(2)}_u(u,1).
\end{split}\label{1.2}
\end{equation}
We have here used a patch numbering as is
indicated in Figure~\ref{fig2}.\vspace{2mm}

Considering a tangent plane continuous 4-patch surface, it must
satisfy that every pair of its patches with a common boundary curve
coincide along that curve. The same must be true of its derivatives
along the same boundary curve. 

Thus

\begin{equation}
\left . \begin{array}{l}
\mathbf{r}^{(1)}_v(1,v) = \mathbf{r}^{(2)}_v(0,v)\\
\mathbf{r}^{(4)}_v(1,v)  = \mathbf{r}^{(3)}_v(0,v),\, v\in[0,1]
\end{array}\right . \label{1.3}
\end{equation}
and
\begin{equation}
\left . \begin{array}{l}
\mathbf{r}^{(1)}_u(u,1)  = \mathbf{r}^{(4)}_u(u,0)\\
\mathbf{r}^{(2)}_u(u,1)  = \mathbf{r}^{(3)}_u (u,0),\, u\in[0,1].
\end{array}\right .
\label{1.4}
\end{equation}
Using the equations (\ref{1.1})--(\ref{1.4}), we prove the next
theorem, which is the main result in this section. In the case where
the patches are described by polynomials, this result has already been
published by B\'ezier in 1986. See \cite{bezier}, p 44-46. Here we
formulate and prove the result for a general parametrization.
\begin{theorem}\label{thm1}
  Let $S$ be a 4-patch surface consisting of $C^{1}_{\#}$-patches
  $\mathbf{r}^{(i)}$, $i=1,\ldots,4$. Let $V$ be the intersection
  point of the four patches.  Then, necessary and sufficient
  conditions in order for the surface $S$ to be tangent plane
  continuous are that there exist continuous functions $\lambda_{ij}$
  and $\kappa_{ij}$, $i=1, j=2,4$ and $i=2,4, j=3$, satisfying the
  equations (\ref{1.1}) and (\ref{1.2}), and the following
  relations
  \begin{equation}
    \begin{array}{l}
        \kappa_{1,2}=\lambda_{1,4}\,\kappa_{4,3}\\
        \kappa_{1,4}=\lambda_{1,2}\,\kappa_{2,3}
      \end{array}
    \label{thm1.1}
  \end{equation}
  and
  \begin{equation}
    \begin{array}{l}
        \lambda_{1,2}-\lambda_{4,3}=\kappa_{1,4}\,\kappa_{4,3}\\
        \lambda_{1,4}-\lambda_{2,3}=\kappa_{1,2}\,\kappa_{2,3}
      \end{array}
    \label{thm1.2}
  \end{equation}
at the vertex $V$.
\end{theorem}

\noindent
\textbf{Remark.} The notations $\kappa_{1,2}$, $\lambda_{1,4}$, $\kappa_{4,3}$, etc, are
to be interpreted as $\kappa_{1,2}(1)$, $\lambda_{1,4}(1)$, $\kappa_{4,3}(0)$, etc.
i.e., as the value of the functions at the vertex $V$.

\paragraph{Proof.} In order to prove the above statement, we must see
under what conditions the equations (\ref{1.1}) and (\ref{1.2}) are all
satisfied at the vertex $V$. We use here the short notation $\mathbf{r}^{(1)}_u$ for $\mathbf{r}^{(1)}_u(u,v)|_V=\mathbf{r}^{(1)}_u(1,1)$, etc.

We start by eliminating $\mathbf{r}^{(2)}_u$, $\mathbf{r}^{(2)}_v$,
$\mathbf{r}^{(4)}_u$ and $\mathbf{r}^{(4)}_v$ in the equations (\ref{1.1}) and
(\ref{1.2}) by using (\ref{1.3}) and (\ref{1.4}) to get
\begin{equation}
  \left . \begin{array}{rrrrl}
      \lambda_{1,2}\mathbf{r}^{(1)}_u &+\, \kappa_{1,2} \mathbf{r}^{(1)}_v& -\,\mathbf{r}^{(3)}_u
      &&= 0\\
      \lambda_{4,3} \mathbf{r}^{(1)}_u &&-\,\mathbf{r}^{(3)}_u &+\, \kappa_{4,3} \mathbf{r}^{(3)}_v&=0\\
      \kappa_{1,4}\mathbf{r}^{(1)}_u &+\,\lambda_{1,4} \mathbf{r}^{(1)}_v &&-\,\mathbf{r}^{(3)}_v &=0\\
      &\lambda_{2,3} \mathbf{r}^{(1)}_v &+\, \kappa_{2,3}\mathbf{r}^{(3)}_u &-\,\mathbf{r}^{(3)}_v  &=0.
    \end{array}\right . \label{1.5}
\end{equation}
With the use of (\ref{1.5}) we replace the second and fourth equation in the above system by eliminating $\mathbf{r}^{(3)}_u$ and $\mathbf{r}^{(3)}_v$ to get
\begin{equation*}
  \left . \begin{array}{l}
%  \[
(\lambda_{4,3}-\lambda_{1,2}+\kappa_{1,4}\kappa_{4,3})\mathbf{r}^{(1)}_u+(\lambda_{1,4}\kappa_{4,3}-\kappa_{1,2})\mathbf{r}^{(1)}_v=0\phantom{.}\\
%\]
%
%and
%
%\[
(\lambda_{1,2}\kappa_{2,3}-\kappa_{1,4})\mathbf{r}^{(1)}_u+(\lambda_{2,3}-\lambda_{1,4}+\kappa_{1,2}\kappa_{2,3})\mathbf{r}^{(1)}_v=0.
   \end{array}\right .
\end{equation*}
%\]
%
Since the vectors $\mathbf{r}^{(1)}_u, \mathbf{r}^{(1)}_v$ span the tangent plane, it follows from the above equations that (\ref{thm1.1}) and (\ref{thm1.2}) must hold.
% %
This concludes the proof.\hfill$\square$\vspace{.8cm}

Let us look at some simple consequences of Theorem~1. Obviously, the
functions $\lambda_{ij}$ are not allowed to be zero if tangential
continuity is to be satisfied. Thus, if e.g.\ $\kappa_{4,3}(0)\ne0$ then
$\kappa_{1,2}(1)\ne0$, which follows from (\ref{thm1.1}). On the other
hand, if $\kappa_{1,4}(1)=0$ then also $\kappa_{2,3}(0)=0$. This
situation is exemplified in Figure~3. In general, it follows from
equation (\ref{thm1.1}) that the pair of $\kappa_{ij}$'s in each equality
must both be zero or non-zero.
\begin{figure}[hbtp]\label{fig3}
  \begin{center}
    \leavevmode
    \includegraphics[width=8cm]{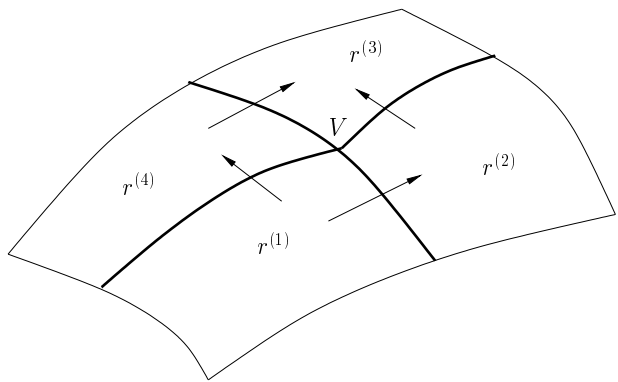}
    \caption{Four patches connected in a common vertex $V$ with
    $\kappa_{1,4}(1)=\kappa_{2,3}(0)=0$ and $\kappa_{1,2}(1)\kappa_{4,3}(0)\ne0$}
  \end{center}
\end{figure}

Another observation that can be done from (\ref{thm1.1}) and
(\ref{thm1.2}) is that the
next relations are true
\begin{equation*}
\begin{split}
\lambda_{1,4}\lambda_{4,3} & =
\lambda_{1,4}(\lambda_{1,2}-\kappa_{1,4}\kappa_{4,3})=\lambda_{1,2}\lambda_{1,4}-\lambda_{1,4}\kappa_{4,3}\kappa_{1,4}\\
& = \lambda_{1,2}\lambda_{1,4}-\kappa_{1,2}\kappa_{1,4}
\end{split}
\end{equation*}
and
\begin{equation*}
\begin{split}
\lambda_{1,2}\lambda_{2,3} & =
\lambda_{1,2}(\lambda_{1,4}-\kappa_{1,2}\kappa_{2,3})=\lambda_{1,2}\lambda_{1,4}-\kappa_{1,2}\lambda_{1,2}\kappa_{2,3}\\
& = \lambda_{1,2}\lambda_{1,4}-\kappa_{1,2}\kappa_{1,4}.
\end{split}
\end{equation*}
In particular, we have
\begin{equation}
\lambda_{1,2}\lambda_{2,3}=\lambda_{1,4}\lambda_{4,3},\label{1.9}
\end{equation}
which will be useful in the next section.

\section{Geometric continuity of order 2}
In this section we consider 4-patch surfaces of higher regularity. Therefore
we introduce the following concept.

\begin{definition}\label{defG2}
  A tangent plane continuous, $G^{1}$, surface is said
  to be curvature continuous, denoted by $G^{2}$, if every point on
  the surface has a unique Dupin indicatrix, which varies continuously
  on the surface. Such a surface is also said to be geometrically
  continuous of order two.
\end{definition}
Another equivalent way to describe the notation of curvature
continuity is to say that the normal curvature at each point and in
each tangential direction\footnote{In fact, it is enough that it holds for 3
  pairwise linearly independent tangential directions by the 3-Tangent Theorem, see Pegna \& Wolter
  \cite{pegnaWolter} or Hoschek \& Lasser \cite{hoschekLasser}, p
  333.} has to be unique, or the principal curvatures are
unique. These differential geometric notations are introduced and
explained in any book about differential geometry, e.g.\ \cite{carmo}.
First we consider a 2-patch surface, where the patches are of regularity
$C^{2}_{\#}=C^{1}_{\#}\cap C^{2}$.
Thus, lack of $G^{2}$-regularity for the surface $S$ can only occur at the common
boundary of the two patches. A necessary and sufficient
condition for a 2-patch surface to be curvature continuous can be found in a paper by Juergen Kahmann \cite{kahmann}. A proof of this well-known result can also be found in \cite{yeLiangNowacki}.
The result is summarized in the next Lemma.
\begin{lemma}\label{lemma2}
  Let $S$ be a 2-patch surface consisting of the $C^{2}_{\#}$-patches
  $(u,v)\mapsto \mathbf{r}^{(1)}(u,v)$ and $(u,v)\mapsto
  \mathbf{r}^{(2)}(u,v)$, $0\le u,v\le1$, satisfying equation
  (\ref{firstEq}) along their common boundary curve. A necessary and
  sufficient condition for the surface $S$ to be curvature continuous
  is that the following relation is fulfilled
  \begin{equation}
    \begin{split}
      \mathbf{r}^{(2)}_{uu}(0,v) & =
      \lambda^2_{1,2}(v)\mathbf{r}^{(1)}_{uu}(1,v)+2\lambda_{1,2}(v)\kappa_{1,2}(v)\mathbf{r}^{(1)}_{uv}(1,v)+\kappa^2_{1,2}(v)\mathbf{r}^{(1)}_{vv}(1,v)\\
      & + \mu_{1,2}(v)\mathbf{r}^{(1)}_u(1,v)+\nu_{1,2}(v)\mathbf{r}^{(1)}_v(1,v),\,\,\,\,
      \mbox{$0\le v\le1$,}
    \end{split}\label{secondEq}
  \end{equation}
where the functions $\lambda_{1,2}$, $\kappa_{1,2}$,
$\mu_{1,2}$ and $\nu_{1,2}$ are continuous.
\end{lemma}

We now continue to consider a 4-patch surface as shown in
Figure~\ref{fig2}. The first and most obvious conditions to have
geometric continuity of order two are that the equations in
(\ref{1.1}) and (\ref{1.2}) are still satisfied when differentiated with
respect to the parameters $u$ and $v$ respectively. We assume here that the functions $\lambda_{ij}$ and $\kappa_{ij}$ are differentiable.
Thus
\begin{equation}
\begin{split}
&
\mathbf{r}^{(2)}_{uv}(0,v)=\lambda'_{1,2}(v)\,\mathbf{r}^{(1)}_{u}(1,v)+\lambda_{1,2}(v)\,\mathbf{r}^{(1)}_{uv}(1,v)+\kappa'_{1,2}(v)\,\mathbf{r}^{(1)}_v(1,v)+\kappa_{1,2}(v)\,\mathbf{r}^{(1)}_{vv}(1,v)\\
&
\mathbf{r}^{(3)}_{uv}(0,v)=\lambda'_{4,3}(v)\,\mathbf{r}^{(4)}_u(1,v)+\lambda_{4,3}(v)\,\mathbf{r}^{(4)}_{uv}(1,v)+\kappa'_{4,3}(v)\,\mathbf{r}^{(4)}_v(1,v)+\kappa_{4,3}(v)\,\mathbf{r}^{(4)}_{vv}(1,v)
\end{split}\label{3.1}
\end{equation}
and
\begin{equation}
\begin{split}
&
\mathbf{r}^{(4)}_{vu}(u,0)=\lambda'_{1,4}(u)\,\mathbf{r}^{(1)}_v(u,1)+\lambda_{1,4}(u)\,\mathbf{r}^{(1)}_{vu}(u,1)+\kappa'_{1,4}(u)\,\mathbf{r}^{(1)}_u(u,1)+\kappa_{1,4}(u)\,\mathbf{r}^{(1)}_{uu}(u,1)\\
&
\mathbf{r}^{(3)}_{vu}(u,0)=\lambda'_{2,3}(u)\,\mathbf{r}^{(2)}_v(u,1)+\lambda_{2,3}(u)\,\mathbf{r}^{(2)}_{vu}(u,1)+\kappa'_{2,3}(u)\,\mathbf{r}^{(2)}_u
(u,1)+\kappa_{2,3}(u)\,\mathbf{r}^{(2)}_{uu}(u,1)
\end{split}\label{3.2}
\end{equation}
must hold. From Lemma~\ref{lemma2}, we know that curvature continuity implies
that the next four relations must be fulfilled

\begin{equation}
\begin{split}
 \mathbf{r}^{(2)}_{uu}(0,v) & =
 \lambda^2_{1,2}(v)\mathbf{r}^{(1)}_{uu}(1,v)+2\lambda_{1,2}(v)\kappa_{1,2}(v)\mathbf{r}^{(1)}_{uv}(1,v)+\kappa^2_{1,2}(v)\mathbf{r}^{(1)}_{vv}(1,v)\\&\makebox[4mm]{\ }+
 \mu_{1,2}(v)\mathbf{r}^{(1)}_u(1,v)+\nu_{1,2}(v)\mathbf{r}^{(1)}_v(1,v)
\end{split}\label{3.3}
\end{equation}
\begin{equation}
\begin{split}
\mathbf{r}^{(3)}_{uu}(0,v) & =
\lambda^2_{4,3}(v)\mathbf{r}^{(4)}_{uu}(1,v)+2\lambda_{4,3}(v)\kappa_{4,3}(v)\mathbf{r}^{(4)}_{uv}(1,v)+\kappa^2_{4,3}(v)\mathbf{r}^{(4)}_{vv}(1,v)\\&\makebox[4mm]{\ }+
 \mu_{4,3}(v)\mathbf{r}^{(4)}_u(1,v)+\nu_{4,3}(v)\mathbf{r}^{(4)}_v(1,v)
\end{split}\label{3.4}
\end{equation}
and
\begin{equation}
\begin{split}
\mathbf{r}^{(3)}_{vv}(u,0) & =
\lambda^2_{2,3}(u)\mathbf{r}^{(2)}_{vv}(u,1)+2\lambda_{2,3}(u)\kappa_{2,3}(u)\mathbf{r}^{(2)}_{uv}(u,1)+\kappa^2_{2,3}(u)\mathbf{r}^{(2)}_{uu}(u,1)\\&\makebox[4mm]{\ }+
 \mu_{2,3}(u)\mathbf{r}^{(2)}_v(u,1)+\nu_{2,3}(u)\mathbf{r}^{(2)}_u(u,1)
\end{split}\label{3.5}
\end{equation}
\begin{equation}
\begin{split}
\mathbf{r}^{(4)}_{vv}(u,0) & =
\lambda^2_{1,4}(u)\mathbf{r}^{(1)}_{vv}(u,1)+2\lambda_{1,4}(u)\kappa_{1,4}(u)\mathbf{r}^{(1)}_{uv}(u,1)+\kappa^2_{1,4}(u)\mathbf{r}^{(1)}_{uu}(u,1)\\&\makebox[4mm]{\ }+
 \mu_{1,4}(u)\mathbf{r}^{(1)}_v(u,1)+\nu_{1,4}(u)\mathbf{r}^{(1)}_u(u,1).
\end{split}\label{3.6}
\end{equation}
The main result of this section is a general result. It concerns parameterized patches of any kind and it is formulated next.
\begin{theorem}
  Let $S$ be a $G^{1}$ 4-patch surface consisting of $C^{2}_{\#}$-patches $\mathbf{r}^{(i)}$, $i=1,\ldots,4$. Let $V$
  be the intersection point of the four patches. Then, necessary and sufficient
  conditions in order for the surface $S$ to be curvature continuous
  %$G^{2}$,
  %%i.e., that the equations in (\ref{3.1})--(\ref{3.6}) are all
  %%satisfied,
  are that there exist continuously differentiable functions $\lambda_{ij}$ and
  $\kappa_{ij}$ and continuous functions $\mu_{ij}$ and
  $\nu_{ij}$, $i=1, j=2,4$ and $i=2,4, j=3$, satisfying the equations (\ref{3.3})--(\ref{3.6})
and that the following relations
\begin{equation}\label{thm2}
\begin{aligned}
&2\lambda_{4,3}\lambda'_{1,4}\kappa_{4,3}-\nu_{1,2}+ \nu_{4,3}\lambda_{1,4}+\mu_{1,4}\kappa^{2}_{4,3}=0\\
&2\lambda_{2,3}\lambda'_{1,2}\kappa_{2,3}-\nu_{1,4}+\nu_{2,3}\lambda_{1,2}+\mu_{1,2}\kappa^{2}_{2,3}=0\\
&2\lambda_{4,3}\kappa_{4,3}\kappa'_{1,4}-\mu_{1,2}+\mu_{4,3}+\nu_{4,3}\kappa_{1,4}+\nu_{1,4}\kappa^{2}_{4,3}=0\\
&2\lambda_{2,3}\kappa_{2,3}\kappa'_{1,2}-\mu_{1,4}+\mu_{2,3}+\nu_{2,3}\kappa_{1,2}+\nu_{1,2}\kappa^{2}_{2,3}=0\\
&\lambda'_{4,3}-\lambda_{2,3}\lambda'_{1,2}+\lambda_{4,3}\kappa'_{1,4}-\lambda_{1,2}\kappa'_{2,3}+\kappa_{1,4}\kappa'_{4,3}-\mu_{1,2}\kappa_{2,3}+\nu_{1,4}\kappa_{4,3}=0\\
&\lambda'_{2,3}-\lambda_{4,3}\lambda'_{1,4}+\lambda_{2,3}\kappa'_{1,2}-\lambda_{1,4}\kappa'_{4,3}+\kappa_{1,2}\kappa'_{2,3}-\mu_{1,4}\kappa_{4,3}+\nu_{1,2}\kappa_{2,3}=0
\end{aligned}
\end{equation}
hold at the vertex $V$.
\end{theorem}
\noindent
\textbf{Remark.} We use here the same short notations as in Theorem~1, i.e., $\lambda_{1,2}$, $\kappa_{1,2}$, $\mu_{1,2}$, $\nu_{1,2}$ etc, are
to be interpreted as $\lambda_{1,2}(1)$, $\kappa_{1,2}(1)$, $\mu_{1,2}(1)$, $\nu_{1,2}(1)$, etc,
i.e., as the value of the functions at the vertex $V$.

\paragraph{Proof}
The idea in the proof is to study the equation system consisting of (\ref{3.1})--(\ref{3.6}) and to reduce, as far as possible, the number of equations including derivatives of second order. To simplify the notations in the proof we use $\mathbf{r}^{(1)}_u$ for $\mathbf{r}^{(1)}_u(u,v)|_V=\mathbf{r}^{(1)}_u(1,1)$, etc.

First, we know that the equations (\ref{3.1}) and (\ref{3.2}) must
hold. Using the relations (\ref{1.3}) and (\ref{1.4}) in order to
replace the vectors $\mathbf{r}^{(4)}_u$ and $\mathbf{r}^{(2)}_v$ with
$\mathbf{r}^{(1)}_u$ and $\mathbf{r}^{(1)}_v$ respectively. Combining this with (\ref{1.1}) and (\ref{1.2}) in order to eliminate $\mathbf{r}^{(2)}_u$ and $\mathbf{r}^{(4)}_v$ we get
\[
\begin{split}
&
\lambda_{1,2}\mathbf{r}^{(1)}_{uv}-\mathbf{r}^{(2)}_{uv}+\lambda'_{1,2}\mathbf{r}^{(1)}_u+\kappa'_{1,2}\mathbf{r}^{(1)}_v+\kappa_{1,2}\mathbf{r}^{(1)}_{vv}=0\phantom{.}\\
&\lambda_{4,3}\mathbf{r}^{(4)}_{uv}-\mathbf{r}^{(3)}_{uv}+\lambda'_{4,3}\mathbf{r}^{(1)}_u+\kappa'_{4,3}(\lambda_{1,4}\mathbf{r}^{(1)}_v+\kappa_{1,4}\mathbf{r}^{(1)}_u)+\kappa_{4,3}\mathbf{r}^{(4)}_{vv}=0\phantom{.}\\
&\lambda_{1,4}\mathbf{r}^{(1)}_{vu}-\mathbf{r}^{(4)}_{vu}+
\kappa'_{1,4}\mathbf{r}^{(1)}_u+\lambda'_{1,4} \mathbf{r}^{(1)}_v+\kappa_{1,4}\mathbf{r}^{(1)}_{uu}=0\phantom{.}\\
&\lambda_{2,3}\mathbf{r}^{(2)}_{vu}-\mathbf{r}^{(3)}_{vu}+\kappa'_{2,3}(\lambda_{1,2}\mathbf{r}^{(1)}_u+\kappa_{1,2}\mathbf{r}^{(1)}_v)+\lambda'_{2,3}\mathbf{r}^{(1)}_v+\kappa_{2,3}\mathbf{r}^{(2)}_{uu}=0.
\end{split}\label{3.7}
\]
In the above equation system we multiply the first equation by $\lambda_{2,3}$, the second one with $-1$, the third one with $-\lambda_{4,3}$ and add them to the fourth equation. We get
%
%\[\left .
\begin{equation}\label{3.21}
 \begin{aligned}
% &\lambda_{1,2}\mathbf{r}^{(1)}_{uv}-\mathbf{r}^{(2)}_{uv}+\lambda'_{1,2}\mathbf{r}^{(1)}_u+\kappa'_{1,2}\mathbf{r}^{(1)}_v+\kappa_{1,2}\mathbf{r}^{(1)}_{vv}=0\\
% &\lambda_{4,3}\mathbf{r}^{(4)}_{uv}-\mathbf{r}^{(3)}_{uv}+\lambda'_{4,3}\mathbf{r}^{(1)}_u+\kappa'_{4,3}\mathbf{r}^{(4)}_v+\kappa_{4,3}\mathbf{r}^{(4)}_{vv}=0\\
% &\lambda_{1,4}\mathbf{r}^{(1)}_{vu}-\mathbf{r}^{(4)}_{vu}+
% \kappa'_{1,4}\mathbf{r}^{(1)}_u+\lambda'_{1,4} \mathbf{r}^{(1)}_v+\kappa_{1,4}\mathbf{r}^{(1)}_{uu}=0\\
(-\lambda'_{4,3}+\lambda_{2,3}\lambda'_{1,2}-\lambda_{4,3}\kappa'_{1,4}+\lambda_{1,2}\kappa'_{2,3}-\kappa_{1,4}\kappa'_{4,3})\mathbf{r}^{(1)}_u+(\lambda'_{2,3}+\lambda_{2,3}\kappa'_{1,2}-\lambda_{4,3}\lambda'_{1,4}\\
\ \ \ \ \ +\kappa_{1,2}\kappa'_{2,3}-\lambda_{1,4}\kappa'_{4,3})\mathbf{r}^{(1)}_v
-\lambda_{4,3}\kappa_{1,4}\mathbf{r}^{(1)}_{uu}+\lambda_{2,3}\kappa_{1,2}\mathbf{r}^{(1)}_{vv}+\kappa_{2,3}\mathbf{r}^{(2)}_{uu}-\kappa_{4,3}\mathbf{r}^{(4)}_{vv}=0.
\end{aligned}%\tag{3.21}
%\right .\]
\end{equation}

Considering the last part in the equation (\ref{3.21}) and
using the equations (\ref{3.3}) and (\ref{3.6}) it follows that
\begin{eqnarray*}
&&\makebox[-15mm]{\ }-\lambda_{4,3}\kappa_{1,4}\mathbf{r}^{(1)}_{uu}+\lambda_{2,3}\kappa_{1,2}\mathbf{r}^{(1)}_{vv}+
\kappa_{2,3}\mathbf{r}^{(2)}_{uu}-\kappa_{4,3}\mathbf{r}^{(4)}_{vv}\\
& =&
-\lambda_{4,3}\kappa_{1,4}\mathbf{r}^{(1)}_{uu}+\lambda_{2,3}\kappa_{1,2}\mathbf{r}^{(1)}_{vv}\\
&& \ \ \ +
\kappa_{2,3}(\lambda^2_{1,2}\mathbf{r}^{(1)}_{uu}+2\lambda_{1,2}\kappa_{1,2}\mathbf{r}^{(1)}_{uv}+
\kappa^2_{1,2}\mathbf{r}^{(1)}_{vv}+\mu_{1,2}\mathbf{r}^{(1)}_{u}+\nu_{1,2}\mathbf{r}^{(1)}_{v})\\
&& \ \ \ -
\kappa_{4,3}(\lambda^2_{1,4}\mathbf{r}^{(1)}_{vv}+2\lambda_{1,4}\kappa_{1,4}\mathbf{r}^{(1)}_{uv}+
\kappa^2_{1,4}\mathbf{r}^{(1)}_{uu}+\mu_{1,4}\mathbf{r}^{(1)}_{v}+\nu_{1,4}\mathbf{r}^{(1)}_{u})\\
&=&
\kappa_{2,3}(\mu_{1,2}\mathbf{r}^{(1)}_u+\nu_{1,2}\mathbf{r}^{(1)}_v)-\kappa_{4,3}(\mu_{1,4}\mathbf{r}^{(1)}_v+\nu_{1,4}\mathbf{r}^{(1)}_u)\\
&& \ \ \ +
(\kappa_{2,3}\lambda^2_{1,2}-\kappa_{4,3}\kappa^2_{1,4}-\lambda_{4,3}\kappa_{1,4})\mathbf{r}^{(1)}_{uu}\\
&& \ \ \ +
(\kappa_{2,3}\kappa^2_{1,2}-\kappa_{4,3}\lambda^2_{1,4}+\lambda_{2,3}\kappa_{1,2})\mathbf{r}^{(1)}_{vv}\\
&& \ \ \ + 2(\kappa_{2,3}\lambda_{1,2}\kappa_{1,2}-\kappa_{4,3}\lambda_{1,4}\kappa_{1,4})\mathbf{r}^{(1)}_{uv}.
\end{eqnarray*}

In order to further reduce the above formula, we use the relations (\ref{thm1.1}) and
(\ref{thm1.2}) to get
\begin{eqnarray*}
\kappa_{2,3}\lambda^2_{1,2}-\kappa_{4,3}\kappa^2_{1,4}-\lambda_{4,3}\kappa_{1,4}&=&\kappa_{1,4}\lambda_{1,2}-\kappa_{4,3}\kappa^2_{1,4}-\lambda_{4,3}\kappa_{1,4}\\
& =& \kappa_{1,4}(\lambda_{1,2}-\lambda_{4,3}-\kappa_{1,4}\kappa_{4,3})=0,\\
%\\
\kappa_{2,3}\kappa^2_{1,2}-\kappa_{4,3}\lambda^2_{1,4}+\lambda_{2,3}\kappa_{1,2}&
=&
\kappa_{2,3}\kappa^2_{1,2}-\kappa_{1,2}\lambda_{1,4}+\lambda_{2,3}\kappa_{1,2}\\
& = &
\kappa_{1,2}(\lambda_{2,3}-\lambda_{1,4}+\kappa_{1,2}\kappa_{2,3})=0\\
\noalign{and}
\kappa_{2,3}\lambda_{1,2}\kappa_{1,2}-\kappa_{4,3}\lambda_{1,4}\kappa_{1,4}&=&\kappa_{1,4}\kappa_{1,2}-\kappa_{1,2}\kappa_{1,4}=0.
\end{eqnarray*}
Thus, it follows that
\[
\begin{split}
&
-\lambda_{4,3}\kappa_{1,4}\mathbf{r}^{(1)}_{uu}+\lambda_{2,3}\kappa_{1,2}\mathbf{r}^{(1)}_{vv}+\kappa_{2,3}\mathbf{r}^{(3)}_{uu}-\kappa_{4,3}\mathbf{r}^{(3)}_{vv}\\
& \phantom{some}=
(\kappa_{2,3}\mu_{1,2}-\kappa_{4,3}\nu_{1,4})\mathbf{r}^{(1)}_{u}+(\kappa_{2,3}\nu_{1,2}-\kappa_{4,3}\mu_{1,4})\mathbf{r}^{(1)}_{v}.
\end{split}
\]
Input the above equality into the equation (\ref{3.21}). As before, the
independence of the two vectors $\mathbf{r}^{(1)}_{u}$ and $\mathbf{r}^{(1)}_{v}$
gives
\begin{equation}
\begin{split}
-\lambda'_{4,3}-\lambda_{4,3}\kappa'_{1,4}+\lambda_{2,3}\lambda'_{1,2}+\lambda_{1,2}\kappa'_{2,3}-\kappa_{1,4}\kappa'_{4,3}+\kappa_{2,3}\mu_{1,2}-\kappa_{4,3}\nu_{1,4}=0\hphantom{.}\\
\lambda'_{2,3}-\lambda_{4,3}\lambda'_{1,4}+\lambda_{2,3}\kappa'_{1,2}+\kappa_{1,2}\kappa'_{2,3}-\lambda_{1,4}\kappa'_{4,3}+\kappa_{2,3}\nu_{1,2}-\kappa_{4,3}\mu_{1,4}=0.
\end{split}\label{3.22}
\end{equation}
The equations in (\ref{3.22}) are necessary in order to fulfill the condition of
geometric continuity of order 2.

Let us next study the equations (\ref{3.3})--(\ref{3.6}) more
closely. We rewrite two of these equations. Start with the equation (\ref{3.4}) by adding to it the first equation in (\ref{3.1}) multiplied by $2\lambda_{4,3}\,\kappa_{4,3}$ and equation (\ref{3.3}) multiplied by $-\lambda_{4,3}/\lambda_{1,2}$. We get the new equation
%where also the equalities in system (\ref{3.21}) have been used, we get the equivalent equation system
\begin{equation}
%\[\left .
\begin{aligned}
%&
%2\lambda_{1,2}\kappa_{1,2}\mathbf{r}^{(1)}_{uv}+\mu_{1,2}\mathbf{r}^{(1)}_{u}+\nu_{1,2}\mathbf{r}^{(1)}_{v}+\lambda^{2}_{1,2}\mathbf{r}^{(1)}_{uu}+\kappa^{2}_{1,2}\mathbf{r}^{(1)}_{vv}-\mathbf{r}^{(3)}_{uu}=0\phantom{.}\\
&
(\mu_{4,3}+2\lambda_{4,3}\kappa_{4,3}\kappa'_{1,4}-\mu_{1,2}\frac{\lambda_{4,3}}{\lambda_{1,2}})\mathbf{r}^{(1)}_{u}+(2\lambda_{4,3}\kappa_{4,3}\lambda'_{1,4}-\nu_{1,2}\frac{\lambda_{4,3}}{\lambda_{1,2}})\mathbf{r}^{(1)}_{v}+\nu_{4,3}\mathbf{r}^{(3)}_{v}\\
&\makebox[2mm]{\ }+(\lambda^{2}_{4,3}+2\lambda_{4,3}\kappa_{4,3}\kappa_{1,4}-\lambda_{1,2}\lambda_{4,3})\mathbf{r}^{(1)}_{uu}-\kappa^{2}_{1,2}\frac{\lambda_{4,3}}{\lambda_{1,2}}\mathbf{r}^{(1)}_{vv}
-(1-\frac{\lambda_{4,3}}{\lambda_{1,2}})\mathbf{r}^{(3)}_{uu}+\kappa^{2}_{4,3}\mathbf{r}^{(3)}_{vv}=0.
%&
%(2\lambda'_{1,2}\lambda_{2,3}\kappa_{2,3}-\nu_{1,4}\frac{\lambda_{2,3}}{\lambda_{1,4}})\mathbf{r}^{(1)}_{u}+(\mu_{2,3}+2\kappa'_{1,2}\lambda_{2,3}\kappa_{2,3}-\mu_{1,4}\frac{\lambda_{2,3}}{\lambda_{1,4}})\mathbf{r}^{(1)}_{v}+\nu_{2,3}\mathbf{r}^{(3)}_{u}\\
%&\makebox[5mm]{\ }-\kappa^{2}_{1,4}\frac{\lambda_{2,3}}{\lambda_{1,4}}\mathbf{r}^{(1)}_{uu}+(\lambda^{2}_{2,3}+2\kappa_{1,2}\lambda_{2,3}\kappa_{2,3}-\lambda_{1,4}\lambda_{2,3})\mathbf{r}^{(1)}_{vv}+\kappa^{2}_{2,3}\mathbf{r}^{(3)}_{uu}-(1-\frac{\lambda_{2,3}}{\lambda_{1,4}})\mathbf{r}^{(3)}_{vv}=0\\
%&
%2\lambda_{1,4}\kappa_{1,4}\mathbf{r}^{(1)}_{uv}+\nu_{1,4}\mathbf{r}^{(1)}_u+\mu_{1,4}\mathbf{r}^{(1)}_v+\kappa^2_{1,4}\mathbf{r}^{(1)}_{uu}+\lambda^2_{1,4}(u)\mathbf{r}^{(1)}_{vv}-\mathbf{r}^{(3)}_{vv}=0.
%% &
%% (\nu_{1,4}-\mu_{1,2}\frac{\lambda_{1,4}\kappa_{1,4}}{\lambda_{1,2}\kappa_{1,2}})\mathbf{r}^{(1)}_{u}+(\mu_{1,4}-\nu_{1,2}\frac{\lambda_{1,4}\kappa_{1,4}}{\lambda_{1,2}\kappa_{1,2}})\mathbf{r}^{(1)}_{v}+(\kappa^{2}_{1,4}-\frac{\lambda_{1,2}\lambda_{1,4}\kappa_{1,4}}{\kappa_{1,2}})\mathbf{r}^{(1)}_{uu}\\
%% &\makebox[5mm]{\ }+(\lambda^{2}_{1,4}-\frac{\lambda_{1,4}\kappa_{1,2}\kappa_{1,4}}{\lambda_{1,2}})\mathbf{r}^{(1)}_{vv}+\frac{\lambda_{1,4}\kappa_{1,4}}{\lambda_{1,2}\kappa_{1,2}}\mathbf{r}^{(3)}_{uu}-\mathbf{r}^{(3)}_{vv}=0.\\
\end{aligned}\label{3.23}
%\right .\]
\end{equation}
We continue with the equation (\ref{3.5}). To this one we add the first equation in (\ref{3.2}) multiplied by $2\lambda_{2,3}\,\kappa_{2,3}$ and equation (\ref{3.6}) multiplied by $-\lambda_{2,3}/\lambda_{1,4}$. We have
\begin{equation}
\begin{aligned}
&
(2\lambda'_{1,2}\lambda_{2,3}\kappa_{2,3}-\nu_{1,4}\frac{\lambda_{2,3}}{\lambda_{1,4}})\mathbf{r}^{(1)}_{u}+(\mu_{2,3}+2\kappa'_{1,2}\lambda_{2,3}\kappa_{2,3}-\mu_{1,4}\frac{\lambda_{2,3}}{\lambda_{1,4}})\mathbf{r}^{(1)}_{v}+\nu_{2,3}\mathbf{r}^{(3)}_{u}\\
&\makebox[2mm]{\ }-\kappa^{2}_{1,4}\frac{\lambda_{2,3}}{\lambda_{1,4}}\mathbf{r}^{(1)}_{uu}+(\lambda^{2}_{2,3}+2\kappa_{1,2}\lambda_{2,3}\kappa_{2,3}-\lambda_{1,4}\lambda_{2,3})\mathbf{r}^{(1)}_{vv}+\kappa^{2}_{2,3}\mathbf{r}^{(3)}_{uu}-(1-\frac{\lambda_{2,3}}{\lambda_{1,4}})\mathbf{r}^{(3)}_{vv}=0.
\end{aligned}\label{3.24}
\end{equation}

Let us now consider the equation (\ref{3.23}). We
want to rewrite this equation in order to make it easier to
handle. Using the relations (\ref{1.2})--(\ref{1.3}), (\ref{thm1.1})--(\ref{thm1.2})  combined with
(\ref{3.3}) and (\ref{3.6}), we get
\begin{equation*}
\begin{split}
&
(\mu_{4,3}\lambda_{1,2}+2\lambda_{1,2}\lambda_{4,3}\kappa_{4,3}\kappa'_{1,4}-\mu_{1,2}\lambda_{4,3})\mathbf{r}^{(1)}_{u}
+\lambda_{4,3}(2\lambda_{1,2}\kappa_{4,3}\lambda'_{1,4}-\nu_{1,2})\mathbf{r}^{(1)}_{v}
+\nu_{4,3}\lambda_{1,2}\mathbf{r}^{(4)}_{v}\\&\makebox[4mm]{\ }
+\lambda_{1,2}\lambda_{4,3}(\lambda_{4,3}+2\kappa_{1,4}\kappa_{4,3}
-\lambda_{1,2})\mathbf{r}^{(1)}_{uu}-\lambda_{4,3}\kappa^{2}_{1,2}\mathbf{r}^{(1)}_{vv}
-(\lambda_{1,2}-\lambda_{4,3})\mathbf{r}^{(2)}_{uu}+\lambda_{1,2}\kappa^{2}_{4,3}\mathbf{r}^{(4)}_{vv}\\[1mm]
&
=(\mu_{4,3}\lambda_{1,2}+2\lambda_{1,2}\lambda_{4,3}\kappa_{4,3}\kappa'_{1,4}
-\mu_{1,2}\lambda_{4,3})\mathbf{r}^{(1)}_{u}
+\lambda_{4,3}(2\lambda_{1,2}\kappa_{4,3}\lambda'_{1,4}-\nu_{1,2})\mathbf{r}^{(1)}_{v}
\\&\makebox[4mm]{\ }+\nu_{4,3}\lambda_{1,2}(\kappa_{1,4}\mathbf{r}^{(1)}_{u}+\lambda_{1,4}\mathbf{r}^{(1)}_{v})
+\lambda_{1,2}\lambda_{4,3}\kappa_{1,4}\kappa_{4,3}\mathbf{r}^{(1)}_{uu}
-\lambda_{4,3}\kappa^{2}_{1,2}\mathbf{r}^{(1)}_{vv}
-\kappa_{1,4}\kappa_{4,3}(\lambda^2_{1,2}\mathbf{r}^{(1)}_{uu}\\&\makebox[4mm]{\ }+2\lambda_{1,2}\kappa_{1,2}\mathbf{r}^{(1)}_{uv}
+\kappa^2_{1,2}\mathbf{r}^{(1)}_{vv}+ \mu_{1,2}\mathbf{r}^{(1)}_u+\nu_{1,2}\mathbf{r}^{(1)}_v)+\lambda_{1,2}\kappa^{2}_{4,3}(\lambda^2_{1,4}\mathbf{r}^{(1)}_{vv}
+2\lambda_{1,4}\kappa_{1,4}\mathbf{r}^{(1)}_{uv}\\&\makebox[4mm]{\ }+\kappa^2_{1,4}\mathbf{r}^{(1)}_{uu}
+ \mu_{1,4}\mathbf{r}^{(1)}_v+\nu_{1,4}\mathbf{r}^{(1)}_u)\\[1mm]
&
=(\mu_{4,3}\lambda_{1,2}+2\lambda_{1,2}\lambda_{4,3}\kappa_{4,3}\kappa'_{1,4}-
\mu_{1,2}\lambda_{4,3}+\nu_{4,3}\lambda_{1,2}\kappa_{1,4}-\kappa_{1,4}\kappa_{4,3}\mu_{1,2}
+\lambda_{1,2}\kappa^{2}_{4,3}\nu_{1,4})\mathbf{r}^{(1)}_{u}
\\&\makebox[4mm]{\ }+(2\lambda_{1,2}\lambda_{4,3}\kappa_{4,3}\lambda'_{1,4}-\nu_{1,2}\lambda_{4,3}+
\nu_{4,3}\lambda_{1,2}\lambda_{1,4}-\kappa_{1,4}\kappa_{4,3}\nu_{1,2}
+\lambda_{1,2}\kappa^{2}_{4,3}\mu_{1,4})\mathbf{r}^{(1)}_{v}
\\&\makebox[4mm]{\ }+\lambda_{1,2}\kappa_{1,4}\kappa_{4,3}(\lambda_{4,3}-\lambda_{1,2}
+\kappa_{1,4}\kappa_{4,3})\mathbf{r}^{(1)}_{uu}-\kappa^{2}_{1,2}(\lambda_{4,3}+\kappa_{1,4}\kappa_{4,3}-\lambda_{1,2})\mathbf{r}^{(1)}_{vv}\\&\makebox[4mm]{\ }
-2\lambda_{1,2}\kappa_{1,4}\kappa_{4,3}(\kappa_{1,2}-\lambda_{1,4}\kappa_{4,3})\mathbf{r}^{(1)}_{uv}\\[1mm]
&
=(\mu_{4,3}\lambda_{1,2}+2\lambda_{1,2}\lambda_{4,3}\kappa_{4,3}\kappa'_{1,4}-
\mu_{1,2}\lambda_{4,3}+\nu_{4,3}\lambda_{1,2}\kappa_{1,4}-\kappa_{1,4}\kappa_{4,3}\mu_{1,2}
+\lambda_{1,2}\kappa^{2}_{4,3}\nu_{1,4})\mathbf{r}^{(1)}_{u} \\&\makebox[4mm]{\ }
+(2\lambda_{1,2}\lambda_{4,3}\kappa_{4,3}\lambda'_{1,4}-\nu_{1,2}\lambda_{4,3}+
\nu_{4,3}\lambda_{1,2}\lambda_{1,4}-\kappa_{1,4}\kappa_{4,3}\nu_{1,2}
+\lambda_{1,2}\kappa^{2}_{4,3}\mu_{1,4})\mathbf{r}^{(1)}_{v}=0.
\end{split}
\end{equation*}
The independence of the tangential vectors $\mathbf{r}^{(1)}_{u}$ and
$\mathbf{r}^{(1)}_{v}$ combined with the first relation in (\ref{thm1.2}) gives
% \begin{equation*}
% \begin{split}
%   \mu_{4,3}\lambda_{1,2}+2\lambda_{1,2}\lambda_{4,3}\kappa_{4,3}\kappa'_{1,4}-
% \mu_{1,2}\lambda_{4,3}+\nu_{4,3}\lambda_{1,2}\kappa_{1,4}-\mu_{1,2}\kappa_{1,4}\kappa_{4,3}
% +\nu_{1,4}\lambda_{1,2}\kappa^{2}_{4,3}=0\\
% 2\lambda_{1,2}\lambda_{4,3}\kappa_{4,3}\lambda'_{1,4}-\nu_{1,2}\lambda_{4,3}+
% \nu_{4,3}\lambda_{1,2}\lambda_{1,4}-\nu_{1,2}\kappa_{1,4}\kappa_{4,3}
% +\mu_{1,4}\lambda_{1,2}\kappa^{2}_{4,3}=0
% \end{split}
% \end{equation*}
% or equivalently, by the use of relation (\ref{thm1.2}), we get
\begin{equation}
\begin{split}
\lambda_{1,2}(\mu_{4,3}+2\lambda_{4,3}\kappa_{4,3}\kappa'_{1,4}-
\mu_{1,2}+\nu_{4,3}\kappa_{1,4}+\nu_{1,4}\kappa^{2}_{4,3})=0\phantom{.}\\
\lambda_{1,2}(2\lambda_{4,3}\kappa_{4,3}\lambda'_{1,4}-\nu_{1,2}+
\nu_{4,3}\lambda_{1,4}+\mu_{1,4}\kappa^{2}_{4,3})=0.
\end{split}\label{3.25}
\end{equation}
We end our examination by simplifying the equation (\ref{3.24}). Similarly as above, we use (\ref{1.1}), (\ref{1.4}),
(\ref{thm1.1})--(\ref{thm1.2}),  combined with
(\ref{3.3}) and (\ref{3.6}). We have
\begin{equation*}
\begin{split}
&
\lambda_{2,3}(2\lambda'_{1,2}\lambda_{1,4}\kappa_{2,3}-\nu_{1,4})\mathbf{r}^{(1)}_{u}
+(\mu_{2,3}\lambda_{1,4}+2\kappa'_{1,2}\lambda_{1,4}\lambda_{2,3}\kappa_{2,3}
-\mu_{1,4}\lambda_{2,3})\mathbf{r}^{(1)}_{v}
+\nu_{2,3}\lambda_{1,4}\mathbf{r}^{(2)}_{u}
\\&\makebox[4mm]{\ }-\lambda_{2,3}\kappa^{2}_{1,4}\mathbf{r}^{(1)}_{uu}+\lambda_{1,4}\lambda_{2,3}(\lambda_{2,3}+2\kappa_{1,2}\kappa_{2,3}-\lambda_{1,4})\mathbf{r}^{(1)}_{vv}
+\lambda_{1,4}\kappa^{2}_{2,3}\mathbf{r}^{(2)}_{uu}-(\lambda_{1,4}-\lambda_{2,3})\mathbf{r}^{(4)}_{vv}\\[1mm]
&
=\lambda_{2,3}(2\lambda'_{1,2}\lambda_{1,4}\kappa_{2,3}-\nu_{1,4})\mathbf{r}^{(1)}_{u}
+(\mu_{2,3}\lambda_{1,4}+2\kappa'_{1,2}\lambda_{1,4}\lambda_{2,3}\kappa_{2,3}
-\mu_{1,4}\lambda_{2,3})\mathbf{r}^{(1)}_{v}
\\&\makebox[4mm]{\ }+\nu_{2,3}\lambda_{1,4}(\lambda_{1,2}\mathbf{r}^{(1)}_{u}+\kappa_{1,2}\mathbf{r}^{(1)}_{v})
-\lambda_{2,3}\kappa^{2}_{1,4}\mathbf{r}^{(1)}_{uu}+
\lambda_{1,4}\lambda_{2,3}\kappa_{1,2}\kappa_{2,3}\mathbf{r}^{(1)}_{vv}
%
%+\lambda_{1,4}\kappa^{2}_{2,3}\mathbf{r}^{(2)}_{uu}\\&\makebox[4mm]{\ }
%-\kappa_{1,2}\kappa_{2,3}\mathbf{r}^{(4)}_{vv}\\[1mm]
+\lambda_{1,4}\kappa^{2}_{2,3}(\lambda^2_{1,2}\mathbf{r}^{(1)}_{uu}
\\&\makebox[4mm]{\ }+2\lambda_{1,2}\kappa_{1,2}\mathbf{r}^{(1)}_{uv}
+\kappa^2_{1,2}\mathbf{r}^{(1)}_{vv}
+ \mu_{1,2}\mathbf{r}^{(1)}_u+\nu_{1,2}\mathbf{r}^{(1)}_v)-\kappa_{1,2}\kappa_{2,3}(\lambda^2_{1,4}\mathbf{r}^{(1)}_{vv}
+2\lambda_{1,4}\kappa_{1,4}\mathbf{r}^{(1)}_{uv}\\&\makebox[4mm]{\ }+\kappa^2_{1,4}\mathbf{r}^{(1)}_{uu}+
\mu_{1,4}\mathbf{r}^{(1)}_v+\nu_{1,4}\mathbf{r}^{(1)}_u)\\[1mm]
&
=(2\lambda_{2,3}\lambda'_{1,2}\lambda_{1,4}\kappa_{2,3}-\nu_{1,4}\lambda_{2,3}
+\nu_{2,3}\lambda_{1,4}\lambda_{1,2}+\mu_{1,2}\lambda_{1,4}\kappa^{2}_{2,3}
-\nu_{1,4}\kappa_{1,2}\kappa_{2,3})\mathbf{r}^{(1)}_{u}
+(\mu_{2,3}\lambda_{1,4}
\\&\makebox[4mm]{\ }+2\kappa'_{1,2}\lambda_{1,4}\lambda_{2,3}\kappa_{2,3}-\mu_{1,4}\lambda_{2,3}+\nu_{2,3}\lambda_{1,4}\kappa_{1,2}
+\nu_{1,2}\lambda_{1,4}\kappa^{2}_{2,3}-\mu_{1,4}\kappa_{1,2}\kappa_{2,3})\mathbf{r}^{(1)}_{v}
+\kappa^{2}_{1,4}(\lambda_{1,4}
\\&\makebox[4mm]{\ }-\lambda_{2,3}-\kappa_{1,2}\kappa_{2,3})\mathbf{r}^{(1)}_{uu}+\lambda_{1,4}\kappa_{1,2}\kappa_{2,3}
(\lambda_{2,3}+\kappa_{1,2}\kappa_{2,3}-\lambda_{1,4})\mathbf{r}^{(1)}_{vv}
+2\lambda_{1,4}\kappa_{1,2}\kappa_{2,3}
(\kappa_{2,3}\lambda_{1,2}\\&\makebox[4mm]{\ }-\kappa_{1,4})\mathbf{r}^{(1)}_{uv}\\[1mm]
&
=(2\lambda_{2,3}\lambda'_{1,2}\lambda_{1,4}\kappa_{2,3}-\nu_{1,4}\lambda_{2,3}
+\nu_{2,3}\lambda_{1,4}\lambda_{1,2}+\mu_{1,2}\lambda_{1,4}\kappa^{2}_{2,3}
-\nu_{1,4}\kappa_{1,2}\kappa_{2,3})\mathbf{r}^{(1)}_{u}
+(\mu_{2,3}\lambda_{1,4}\\&\makebox[4mm]{\ }
+2\kappa'_{1,2}\lambda_{1,4}\lambda_{2,3}\kappa_{2,3}-\mu_{1,4}\lambda_{2,3}+\nu_{2,3}\lambda_{1,4}\kappa_{1,2}
+\nu_{1,2}\lambda_{1,4}\kappa^{2}_{2,3}-\mu_{1,4}\kappa_{1,2}\kappa_{2,3})\mathbf{r}^{(1)}_{v}=0.
\end{split}
\end{equation*}
Using the same argument as before, i.e., the independence of the tangential
vectors $\mathbf{r}^{(1)}_{u}$ and $\mathbf{r}^{(1)}_{v}$, together with the second relation in (\ref{thm1.2}), it implies
\begin{equation}
\begin{split}
\lambda_{1,4}(2\lambda_{2,3}\lambda'_{1,2}\kappa_{2,3}-\nu_{1,4}
+\nu_{2,3}\lambda_{1,2}+\mu_{1,2}\kappa^{2}_{2,3})=0\phantom{.}\\
\lambda_{1,4}(\mu_{2,3}+2\kappa'_{1,2}\lambda_{2,3}\kappa_{2,3}
-\mu_{1,4}+\nu_{2,3}\kappa_{1,2}+\nu_{1,2}\kappa^{2}_{2,3})=0.
\end{split}\label{3.26}
\end{equation}

Combining the results in (\ref{3.22}), (\ref{3.25}) and (\ref{3.26})
with the fact that $\lambda_{ij}\ne0$, we get the compatibility
conditions (\ref{thm2}) for $G^{2}$, which are necessary and
sufficient for having a simultaneous satisfaction of the equations
(\ref{3.1})--(\ref{3.6}). This ends the proof.\hfill$\square$

\section{An algorithm}
In the previous part of this article we have achieved general results about regularity of 4-patch surfaces. In this section we will apply those results in the $G^1$-case to a particular patch parametrization such as Bezier representation. Our goal is to make the connection between two patches more flexible, which will make it easier to handle geometrical and other demands. The way to handle this will be done by letting the functions $\lambda$ and $\kappa$ in (\ref{firstEq}) be suitable polynomials.
We start this section by studying the consequences for the relations between the control points from a 2-patch surface in such a case. This result will then be used in creating a smooth 4-patch surface from an uncomplete such surface consisting of three patches. Here the compatibility conditions in Theorem~\ref{thm1} are of fundamental importance. The next step is to solve the problem of filling a hole in a surface in a smooth way, more precisely, we consider an uncomplete regular 9-patch surface as in Figure~5. %%%%%\ref{fig4}.
To create the interior patch we use the same technique as in creating a 4-patch surface. Finally, what we have done so far will be used in creating a fillet surface, see Figure~6, %\ref{fig5},
which together with its surrounding becomes a smooth surface.

\vspace{1mm}

\begin{figure}[hbtp]\label{fig3.5}
  \begin{center}
    \leavevmode
     \includegraphics[width=3cm]{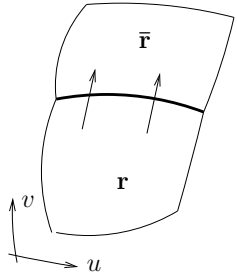}
    \caption{Two patches connected along a common boundary}
  \end{center}
\end{figure}
We begin by considering the well-known $G^1$-relation (\ref{firstEq}) between two patches, see Figure~\ref{fig3.5}. Let the two patches be given by $(u,v)\mapsto \mathbf{r}(u,v)$ and $(u,v)\mapsto \mathbf{\bar r}(u,v)$, where $u,v\in[0,1]$ and with the common boundary curve $u\mapsto \mathbf{\bar r}(u,0)=\mathbf{r}(u,1)$ for $u\in[0,1]$. Then
\begin{equation}\label{eqn4.01}
  \mathbf{\bar r}_v(u,0) =
\lambda(u)\,\mathbf{r}_v(u,1)+\kappa(u)\,\mathbf{r}_u(u,1),\,\,\,\,0\le
u\le1.
\end{equation}
We restrict ourselves in using Bezier representation of the patches as follows
\begin{eqnarray*}
  &\mathbf{r}(u,v) = \sum_{i=0}^{3}\sum_{j=0}^{3} \mathbf{q}_{ij}\,B_i^{(3)}(u)\,B_j^{(3)}(v)\\
  \noalign{and}
  &\mathbf{\bar r}(u,v) = \sum_{i=0}^{5}\sum_{j=0}^{m} \mathbf{\bar q}_{ij}\,B_i^{(5)}(u)\,B_j^{(m)}(v),
\end{eqnarray*}
where the parameters $u,v\in[0,1]$ with $m=4,5,6$. The functions $B_i^{(n)}$ for $i=0,1,\ldots,n$, are the Bernstein polynomials of degree $n$, i.e., $B_i^{(n)}(u)=$ $n\choose i$ $(1-u)^{n-i}\,u^i$ for $u\in[0,1]$.

Let the functions $\lambda$ and $\kappa$ be polynomials satisfying
\begin{equation}\label{eqn4.02}
 \begin{split}
 &\lambda(u)=\lambda_{0}\,(1-u)^2+2\alpha\,(1-u)u+\lambda_1\,u^2\\
%  \noalign{and}&\\
  &\kappa(u)=\kappa_{0}\,(1-u)^3+3\beta^{(1)}\,(1-u)^2u+3\beta^{(2)}\,(1-u)u^2+\kappa_1\,u^3,
\end{split}
\end{equation}
where $\lambda_{0}$, $\lambda_1$, $\kappa_{o}$, $\kappa_1$ $\alpha$, $\beta^{(1)}$ and $\beta^{(2)}$ are constants and the parameter $u\in [0,1]$. From formula (\ref{eqn4.01}) we get
\begin{equation*}
  \begin{split}
    &m \sum_{i=0}^{5} (\mathbf{\bar q}_{i,1}-\mathbf{\bar q}_{i,0})\,B_i^{(5)}(u)
    = (\lambda_0(1-u)^2+2\alpha(1-u)u+\lambda_1 u^2)\sum_{i=0}^{3} 3\,(\mathbf{q}_{i,3}-\mathbf{q}_{i,2})\,B_i^{(3)}(u)
    \\&\makebox[3mm]{}
    +(\kappa_0(1-u)^3+3\beta^{(1)}(1-u)^2u+3\beta^{(2)}(1-u)u^2+\kappa_1\,u^3)\sum_{i=0}^{2}3\, (\mathbf{q}_{i+1,3}-\mathbf{q}_{i,3})\,B_i^{(2)}(u)
    \\&
    = 3\sum_{i=0}^{5} \big(\lambda_0\,(1-\frac{i}{5})(1-\frac{i}{4})\,(\mathbf{q}_{i,3}-\mathbf{q}_{i,2})+2\alpha\,(1-\frac{i}{5})\frac{i}{4}\,(\mathbf{q}_{i-1,3}-\mathbf{q}_{i-1,2})\\&\makebox[3mm]{}
    +\lambda_1\,\frac{i}{4}\frac{i-1}{5}\,(\mathbf{q}_{i-2,3}-\mathbf{q}_{i-2,2})\big)\,B_i^{(5)}(u)
    \\&\makebox[3mm]{}
    +3 \sum_{i=0}^{5} \big(\kappa_0\,(1-\frac{i}{5})(1-\frac{i}{4})(1-\frac{i}{3})\,(\mathbf{q}_{i+1,3}-\mathbf{q}_{i,3})+3\beta^{(1)}\,(1-\frac{i}{5})(1-\frac{i}{4})\frac{i}{3}\,(\mathbf{q}_{i,3}-\mathbf{q}_{i-1,3})
    \\&\makebox[3mm]{}
    +3\beta^{(2)}\,(1-\frac{i}{5})\frac{i}{4}\frac{i-1}{3}\,(\mathbf{q}_{i-1,3}-\mathbf{q}_{i-2,3})+\kappa_1\,\frac{i}{5}\frac{i-1}{4}\frac{i-2}{3}\,(\mathbf{q}_{i-2,3}-\mathbf{q}_{i-3,3})\big)\,B_i^{(5)}(u).
  \end{split}
\end{equation*}
Using the above identity we see that the control points must satisfy the next relations. We have
\begin{equation}
  \begin{split}
    m (\mathbf{\bar q}_{0,1}-\mathbf{\bar q}_{0,0})=& 3\,\big(\lambda_0\,(\mathbf{q}_{0,3}-\mathbf{q}_{0,2})+\kappa_0\,(\mathbf{q}_{1,3}-\mathbf{q}_{0,3})\big)
    \\
    m (\mathbf{\bar q}_{1,1}-\mathbf{\bar q}_{1,0})=&3\,\big(\frac{3\,\lambda_0}{5}\,(\mathbf{q}_{1,3}-\mathbf{q}_{1,2})+\frac{2\alpha}{5}\,(\mathbf{q}_{0,3}-\mathbf{q}_{0,2})\\&
    + \frac{2\,\kappa_0}{5}\,(\mathbf{q}_{2,3}-\mathbf{q}_{1,3})+\frac{3\,\beta^{(1)}}{5}(\mathbf{q}_{1,3}-\mathbf{q}_{0,3})\big)
    \\
    m (\mathbf{\bar q}_{2,1}-\mathbf{\bar q}_{2,0})=&3\,\big(\frac{3\,\lambda_0}{10}\,(\mathbf{q}_{2,3}-\mathbf{q}_{2,2})+\frac{6\,\alpha}{10}\,(\mathbf{q}_{1,3}-\mathbf{q}_{1,2})+\frac{\lambda_1}{10}\,(\mathbf{q}_{0,3}-\mathbf{q}_{0,2})\\&
    +\frac{\kappa_0}{10}\,(\mathbf{q}_{3,3}-\mathbf{q}_{2,3})+\frac{6\,\beta^{(1)}}{10}\,(\mathbf{q}_{2,3}-\mathbf{q}_{1,3})+\frac{3\,\beta^{(2)}}{10}\,(\mathbf{q}_{1,3}-\mathbf{q}_{0,3})\big)
    \\
    m (\mathbf{\bar q}_{3,1}-\mathbf{\bar q}_{3,0})=& 3\,\big(\frac{\lambda_0}{10}\,(\mathbf{q}_{3,3}-\mathbf{q}_{3,2})+\frac{6\,\alpha}{10}\,(\mathbf{q}_{2,3}-\mathbf{q}_{2,2})+\frac{3\,\lambda_1}{10}\,(\mathbf{q}_{1,3}-\mathbf{q}_{1,2})\\&
    +\frac{3\,\beta^{(1)}}{10}\,(\mathbf{q}_{3,3}-\mathbf{q}_{2,3})+\frac{6\,\beta^{(2)}}{10}\,(\mathbf{q}_{2,3}-\mathbf{q}_{1,3})+\frac{\kappa_1}{10}\,(\mathbf{q}_{1,3}-\mathbf{q}_{0,3})\big)
    \\
    m (\mathbf{\bar q}_{4,1}-\mathbf{\bar q}_{4,0})=& 3\,\big(\frac{2\,\alpha}{5}\,(\mathbf{q}_{3,3}-\mathbf{q}_{3,2})+\frac{3\,\lambda_1}{5}\,(\mathbf{q}_{2,3}-\mathbf{q}_{2,2})\\&
    +\frac{3\,\beta^{(2)}}{5}\,(\mathbf{q}_{3,3}-\mathbf{q}_{2,3})+\frac{2\,\kappa_1}{5}\,(\mathbf{q}_{2,3}-\mathbf{q}_{1,3})\big)
    \\
    m (\mathbf{\bar q}_{5,1}-\mathbf{\bar q}_{5,0})=&3\,\big(\lambda_1\,(\mathbf{q}_{3,3}-\mathbf{q}_{3,2})+\kappa_1\,(\mathbf{q}_{3,3}-\mathbf{q}_{2,3})\big).
\end{split}\label{4.05}
\end{equation}

Since the two patches $\mathbf{r}$ and $\mathbf{\bar r}$ coincide along their common boundary curve we also have the following identity
\begin{equation*}
\begin{split}
  &\sum_{i=0}^{5} \mathbf{\bar q}_{i,0}\,B_i^{(5)}(u)= \sum_{i=0}^{3} \mathbf{q}_{i,3}\,B_i^{(3)}(u)\\&=\sum_{i=0}^{5} \big((1-\frac i5)(1-\frac i4)\,\mathbf{q}_{i,3}+2(1-\frac i5)\frac i4\,\mathbf{q}_{i-1,3}+\frac i5\frac{i-1}4\,\mathbf{q}_{i-2,3}\big)\,B_i^{(5)}(u),\,\,\,u\in[0,1],
\end{split}
\end{equation*}
which, by identification, implies that
\begin{equation}
  \begin{split}
    \mathbf{\bar q}_{0,0}&= \mathbf{q}_{0,3}\\
    \mathbf{\bar q}_{1,0}&= \frac35\,\mathbf{q}_{1,3}+\frac 25\,\mathbf{q}_{0,3}\\%\nonumber
    \mathbf{\bar q}_{2,0}&= \frac 3{10}\,\mathbf{q}_{2,3}+\frac 6{10}\,\mathbf{q}_{1,3}+\frac 1{10}\,\mathbf{q}_{0,3}\\
    \mathbf{\bar q}_{3,0}&= \frac 1{10}\,\mathbf{q}_{3,3}+\frac 6{10}\,\mathbf{q}_{2,3}+\frac 3{10}\,\mathbf{q}_{1,3}\\
    \mathbf{\bar q}_{4,0}&= \frac 25\mathbf{q}_{3,3}+\frac 35\mathbf{q}_{2,3}\\
    \mathbf{\bar q}_{5,0}&= \mathbf{q}_{3,3}.
  \end{split}\label{4.06}
\end{equation}
Thus, from the $G^1$-condition between the two patches $\mathbf{r}$ and $\mathbf{\bar r}$, where $\mathbf{r}$ is known, the control points $\mathbf{\bar q}_{i,j}$,\, $i=0,1,\ldots,5, j=0,1$, are forced to satisfy the relations (\ref{4.05}) and (\ref{4.06}). Nevertheless, first we have to decide the value of the parameters $\lambda_{0}$, $\lambda_1$, $\kappa_{0}$, $\kappa_{1}$, $\alpha$, $\beta^{(1)}$ and $\beta^{(2)}$.
In the later applications it will be obvious how to choose certain of these parameters.

In order to define the patch $\mathbf{\bar r}$ completely, the other control points $\mathbf{\bar q}_{i,j}$, $i=0,1,\ldots,5,\, j=2,\ldots,m$, must also be defined but not from the above regularity condition. There must be any other way to define them.

The result we have achieved so far in this section will be used next, where we consider 4-patch surfaces.

\subsection{A 4-patch surface}
In the introduction of this section we considered a 2-patch surface. There we achieved relations between the control points of the two patches in (\ref{4.05}) and (\ref{4.06}). We will use this when creating a 4-patch surface. Looking at Figure~\ref{fig2}, we assume that we have a $G^1$-surface consisting of the three Bezier patches $\mathbf{r}^{(4)}$, $\mathbf{r}^{(1)}$ and $\mathbf{r}^{(2)}$. Our goal is to create the fourth patch $\mathbf{r}^{(3)}$ in such a way that we get a $G^1$-regular 4-patch surface.

Since the patch $\mathbf{r}^{(3)}$ must have a common boundary curve with $\mathbf{r}^{(2)}$ and $\mathbf{r}^{(4)}$ respectively, the second equation in (\ref{1.1}) and (\ref{1.2}) respectively must be fulfilled. We have
\begin{equation}
\begin{split}
\mathbf{r}^{(3)}_v(u,0) & = \lambda_{2,3}(u)\mathbf{r}^{(2)}_v(u,1)+ \kappa_{2,3}(u)\mathbf{r}^{(2)}_u(u,1)\\
\mathbf{r}^{(3)}_u (0,v) & = \lambda_{4,3}(v)\mathbf{r}^{(4)}_u
  (1,v)+\kappa_{4,3}(v)\mathbf{r}^{(4)}_v(1,v),
\end{split}\label{4.11}
\end{equation}
where the parameters $u,v\in[0,1]$.
We continue with the next two equations which we get from the the second equation in (\ref{1.3}) and (\ref{1.4}) respectively. Those equations are the boundary conditions, i.e.,
\begin{equation}
\left . \begin{array}{l}
\mathbf{r}^{(3)}(u,0) = \mathbf{r}^{(2)}(u,1)\\
\mathbf{r}^{(3)}(0,v) = \mathbf{r}^{(4)}(1,v).
\end{array}\right . \label{4.12}
\end{equation}
To be more specific, we exemplify by letting the patches $\mathbf{r}^{(1)}$, $\mathbf{r}^{(2)}$ and $\mathbf{r}^{(4)}$ all be Bezier patches of bi-degree $(3,3)$, while the patch $\mathbf{r}^{(3)}$ is of bi-degree (5,5), i.e.,

\begin{equation*}\label{4.13}
\mathbf{r}^{(3)}(u,v) = \sum_{i=0}^{5}\sum_{j=0}^{5} \mathbf{q}_{ij}^{(3)}\,B_i^{(5)}(u)\,B_j^{(5)}(v),
\end{equation*}
where $u,v\in[0,1]$.

Since the the patches $\mathbf{r}^{(1)}$, $\mathbf{r}^{(2)}$ and $\mathbf{r}^{(4)}$ have the same bi-degree the functions $\lambda_{1,2}$, $\lambda_{1,4}$, $\kappa_{1,2}$ and $\kappa_{1,4}$, defined in the first equation of (\ref{1.1}) and (\ref{1.2}) respectively, must satisfy that $\lambda_{1,2}$ and $\lambda_{1,4}$ are identically constant and $\kappa_{1,2}$ and $\kappa_{1,4}$ are polynomials of first degree. In fact these polynomials must satisfy $\kappa_{1,2}(v)=\kappa_{1,2}(0)(1-v)$ and $\kappa_{1,4}(u)=\kappa_{1,4}(0)(1-u)$, but here we consider the case where the functions $\kappa_{1,2}$ and $\kappa_{1,4}$ are identically zero. There is in principle no difference between the two cases.

Combining the above with the compatibility conditions (\ref{thm1.1}) and (\ref{thm1.2}) we get
\begin{equation}
    \begin{array}{l}
        \lambda_{2,3}(0)=\lambda_{1,4}\\
        \lambda_{4,3}(0)=\lambda_{1,2}\\
        \kappa_{2,3}(0)=0\\
        \kappa_{4,3}(0)=0.
    \end{array} \label{4.14}
\end{equation}

Let the functions $\lambda_{2,3}$, $\lambda_{4,3}$, $\kappa_{2,3}$ and $\kappa_{4,3}$ be defined as in (\ref{eqn4.02}). Using (\ref{4.14}) we have
\begin{equation}\label{4.141}
    \begin{array}{l}
    \lambda_{2,3}(u)=\lambda_{1,4}\,(1-u)^2+2\,\alpha_{2,3}\,(1-u)u+\lambda_{2,3}(1)\,u^2\\
    \lambda_{4,3}(v)=\lambda_{1,2}\,(1-v)^2+2\,\alpha_{4,3}\,(1-v)v+\lambda_{4,3}(1)\,v^2\\
     \kappa_{2,3}(u)=3\,\beta_{2,3}^{(1)}\,(1-u)^2u+3\,\beta_{2,3}^{(2)}\,(1-u)u^2+\kappa_{2,3}(1)\,u^3\\
    %\noalign{and}
    \kappa_{4,3}(v)=3\,\beta_{4,3}^{(1)}\,(1-v)^2v+3\,\beta_{4,3}^{(2)}\,(1-v)v^2+\kappa_{4,3}(1)\,v^3,
    \end{array}
\end{equation}
where $\lambda_{i,3}(1)$, $\kappa_{i,3}(1)$ and $\beta_{i,3}^{(2)}$ for $i=2,4,$ are any constants, while the coefficients $\alpha_{i,3}$ and $\beta_{i,3}^{(1)}$ for $i=2,4$ need to fulfill certain compatibility conditions given in the end of this subsection, see (\ref{4.19}). %We also refer to B\'ezier \cite{bezier} and Sarraga \cite{sarraga}.

The solutions of the system of differential equations (\ref{4.11}) with boundary conditions (\ref{4.12}) together with the compatibility conditions (\ref{4.14}) include the set of all possible patches that will result in a $G^1$-regular $4$-patch surface. In order to get such a solution $\mathbf{r}^{(3)}$ we start by solving the first equation in (\ref{4.12}). From (\ref{4.06}) it follows that
\begin{equation}
  \begin{split}
  \mathbf{q}^{(3)}_{0,0}&= \mathbf{q}^{(2)}_{0,3}\\
  \mathbf{q}^{(3)}_{1,0}&= \frac35\,\mathbf{q}^{(2)}_{1,3}+\frac 25\,\mathbf{q}^{(2)}_{0,3}\\%\nonumber
  \mathbf{q}^{(3)}_{2,0}&= \frac 3{10}\,\mathbf{q}^{(2)}_{2,3}+\frac 6{10}\,\mathbf{q}^{(2)}_{1,3}+\frac 1{10}\,\mathbf{q}^{(2)}_{0,3}\\
  \mathbf{q}^{(3)}_{3,0}&= \frac 1{10}\,\mathbf{q}^{(2)}_{3,3}+\frac 6{10}\,\mathbf{q}^{(2)}_{2,3}+\frac 3{10}\mathbf{q}^{(2)}_{1,3}\\
  \mathbf{q}^{(3)}_{4,0}&= \frac 2{5}\,\mathbf{q}^{(2)}_{3,3}+\frac 3{5}\,\mathbf{q}^{(2)}_{2,3}\\
  \mathbf{q}^{(3)}_{5,0}&= \mathbf{q}^{(2)}_{3,3}.
\end{split}\label{4.15}
\end{equation}
\vspace{2mm}

Considering next the first equation in (\ref{4.11}) combined with (\ref{4.14}), it follows from (\ref{4.05}) in the previous subsection that

\begin{equation}\label{4.16}
  \begin{split}
    \mathbf{q}^{(3)}_{0,1}-\mathbf{q}^{(3)}_{0,0}=& \frac{3\,\lambda_{1,4}}5\,(\mathbf{q}^{(2)}_{0,3}-\mathbf{q}^{(2)}_{0,2})
    \\
    \mathbf{q}^{(3)}_{1,1}-\mathbf{q}^{(3)}_{1,0}=&\frac35\,\big(\frac{3\,\lambda_{1,4}}{5}\,(\mathbf{q}^{(2)}_{1,3}-\mathbf{q}^{(2)}_{1,2})+\frac{2\alpha_{2,3}}{5}\,(\mathbf{q}^{(2)}_{0,3}-\mathbf{q}^{(2)}_{0,2})
    + \frac{3\,\beta_{2,3}^{(1)}}{5}\,(\mathbf{q}^{(2)}_{1,3}-\mathbf{q}^{(2)}_{0,3})\big)
    \\
    \mathbf{q}^{(3)}_{2,1}-\mathbf{q}^{(3)}_{2,0}=&\frac35\,\big(\frac{3\,\lambda_{1,4}}{10}\,(\mathbf{q}^{(2)}_{2,3}-\mathbf{q}^{(2)}_{2,2})+\frac{6\,\alpha_{2,3}}{10}\,(\mathbf{q}^{(2)}_{1,3}-\mathbf{q}^{(2)}_{1,2})+\frac{\lambda_{2,3}(1)}{10}\,(\mathbf{q}^{(2)}_{0,3}-\mathbf{q}^{(2)}_{0,2})\\&
    +\frac{6\,\beta_{2,3}^{(1)}}{10}\,(\mathbf{q}^{(2)}_{2,3}-\mathbf{q}^{(2)}_{1,3})+
    \frac{3\,\beta_{2,3}^{(2)}}{10}\,(\mathbf{q}^{(2)}_{1,3}-\mathbf{q}^{(2)}_{0,3})\big)
    \\
    \mathbf{q}^{(3)}_{3,1}-\mathbf{q}^{(3)}_{3,0}=& \frac35\,\big(\frac{\lambda_{1,4}}{10}\,(\mathbf{q}^{(2)}_{3,3}-\mathbf{q}^{(2)}_{3,2})
    +\frac{6\,\alpha_{2,3}}{10}\,(\mathbf{q}^{(2)}_{2,3}-\mathbf{q}^{(2)}_{2,2})
    +\frac{3\,\lambda_{2,3}(1)}{10}\,(\mathbf{q}^{(2)}_{1,3}-\mathbf{q}^{(2)}_{1,2})\\&
    +\frac{3\,\beta_{2,3}^{(1)}}{10}\,(\mathbf{q}^{(2)}_{3,3}-\mathbf{q}^{(2)}_{2,3})+
    \frac{6\,\beta_{2,3}^{(2)}}{10}\,(\mathbf{q}^{(2)}_{2,3}-\mathbf{q}^{(2)}_{1,3})
    +\frac{\kappa_{2,3}(1)}{10}\,(\mathbf{q}^{(2)}_{1,3}-\mathbf{q}^{(2)}_{0,3})\big)
    \\
    \mathbf{q}^{(3)}_{4,1}-\mathbf{q}^{(3)}_{4,0}=&\frac35\,\big(\frac{2\,\alpha_{2,3}}{5}\,(\mathbf{q}^{(2)}_{3,3}-\mathbf{q}^{(2)}_{3,2})+\frac{3\,\lambda_{2,3}(1)}{5}\,(\mathbf{q}^{(2)}_{2,3}-\mathbf{q}^{(2)}_{2,2})\\&
+\frac{3\,\beta_{2,3}^{(2)}}{5}\,(\mathbf{q}^{(2)}_{3,3}-\mathbf{q}^{(2)}_{2,3})+\frac{2\,\kappa_{2,3}(1)}{5}(\mathbf{q}^{(2)}_{2,3}-\mathbf{q}^{(2)}_{1,3})\big)
    \\
    \mathbf{q}^{(3)}_{5,1}-\mathbf{q}^{(3)}_{5,0}=&\frac35\,\big(\lambda_{2,3}(1)\,(\mathbf{q}^{(2)}_{3,3}-\mathbf{q}^{(2)}_{3,2})+\kappa_{2,3}(1)\,(\mathbf{q}^{(2)}_{3,3}-\mathbf{q}^{(2)}_{2,3})\big).
\end{split}
\end{equation}
When considering the second equation in (\ref{4.11}), using symmetry in (\ref{4.16}) we imme\-dia\-tely get the following relations
\begin{equation}\label{4.17}
  \begin{split}
    \mathbf{q}^{(3)}_{1,0}-\mathbf{q}^{(3)}_{0,0}=& \frac{3\,\lambda_{1,2}}5\,(\mathbf{q}^{(4)}_{3,0}-\mathbf{q}^{(4)}_{2,0})
    \\
    \mathbf{q}^{(3)}_{1,1}-\mathbf{q}^{(3)}_{0,1}=&\frac35\,\big(\frac{3\,\lambda_{1,2}}{5}\,(\mathbf{q}^{(4)}_{3,1}-\mathbf{q}^{(4)}_{2,1})+\frac{2\alpha_{4,3}}{5}\,(\mathbf{q}^{(4)}_{3,0}-\mathbf{q}^{(4)}_{2,0})
    + \frac{3\,\beta_{4,3}^{(1)}}{5}\,(\mathbf{q}^{(4)}_{3,1}-\mathbf{q}^{(4)}_{3,0})\big)
    \\
    \mathbf{q}^{(3)}_{1,2}-\mathbf{q}^{(3)}_{0,2}=&\frac35\,\big(\frac{3\,\lambda_{1,2}}{10}\,(\mathbf{q}^{(4)}_{3,2}-\mathbf{q}^{(4)}_{2,2})+\frac{6\,\alpha_{4,3}}{10}\,(\mathbf{q}^{(4)}_{3,1}-\mathbf{q}^{(4)}_{2,1})+\frac{\lambda_{4,3}(1)}{10}\,(\mathbf{q}^{(4)}_{3,0}-\mathbf{q}^{(4)}_{2,0})\\&
         +\frac{6\,\beta_{4,3}^{(1)}}{10}\,(\mathbf{q}^{(4)}_{3,2}-\mathbf{q}^{(4)}_{3,1})+
        \frac{3\,\beta_{4,3}^{(2)}}{10}\,(\mathbf{q}^{(4)}_{3,1}-\mathbf{q}^{(4)}_{3,0})\big)
    \\
    \mathbf{q}^{(3)}_{1,3}-\mathbf{q}^{(3)}_{0,3}=& \frac35\,\big(\frac{\lambda_{1,2}}{10}\,(\mathbf{q}^{(4)}_{3,3}-\mathbf{q}^{(4)}_{2,3})
    +\frac{6\,\alpha_{4,3}}{10}\,(\mathbf{q}^{(4)}_{3,2}-\mathbf{q}^{(4)}_{2,2})
    +\frac{3\,\lambda_{4,3}(1)}{10}\,(\mathbf{q}^{(4)}_{3,1}-\mathbf{q}^{(4)}_{2,1})\\&
    +\frac{3\,\beta_{4,3}^{(1)}}{10}\,(\mathbf{q}^{(4)}_{3,3}-\mathbf{q}^{(4)}_{3,2})+
    \frac{6\,\beta_{4,3}^{(2)}}{10}\,(\mathbf{q}^{(4)}_{3,2}-\mathbf{q}^{(4)}_{3,1})
    +\frac{\kappa_{4,3}(1)}{10}\,(\mathbf{q}^{(4)}_{3,1}-\mathbf{q}^{(4)}_{3,0})\big)
    \\
    \mathbf{q}^{(3)}_{1,4}-\mathbf{q}^{(3)}_{0,4}=&\frac35\,\big(\frac{2\,\alpha_{4,3}}{5}\,(\mathbf{q}^{(4)}_{3,3}-\mathbf{q}^{(4)}_{2,3})+\frac{3\,\lambda_{4,3}(1)}{5}\,(\mathbf{q}^{(4)}_{3,2}-\mathbf{q}^{(4)}_{2,2})\\&
+\frac{3\,\beta_{4,3}^{(2)}}{5}\,(\mathbf{q}^{(4)}_{3,3}-\mathbf{q}^{(4)}_{3,2})+\frac{2\,\kappa_{4,3}(1)}{5}(\mathbf{q}^{(4)}_{3,2}-\mathbf{q}^{(4)}_{3,1})\big)
    \\
    \mathbf{q}^{(3)}_{1,5}-\mathbf{q}^{(3)}_{0,5}=&\frac35\,\big(\lambda_{4,3}(1)\,(\mathbf{q}^{(4)}_{3,3}-\mathbf{q}^{(4)}_{2,3})+\kappa_{4,3}(1)\,(\mathbf{q}^{(4)}_{3,3}-\mathbf{q}^{(4)}_{3,2})\big).
\end{split}
\end{equation}

Finally, from the second equation in (\ref{4.12}), similarly as in (\ref{4.15}) combined with symmetry, we get
\begin{equation}\label{4.18}
  \begin{split}
  \mathbf{q}^{(3)}_{0,0}&= \mathbf{q}^{(4)}_{3,0}\\
  \mathbf{q}^{(3)}_{0,1}&= \frac35\,\mathbf{q}^{(4)}_{3,1}+\frac 25\,\mathbf{q}^{(4)}_{3,0}\\%\nonumber
  \mathbf{q}^{(3)}_{0,2}&= \frac 3{10}\,\mathbf{q}^{(4)}_{3,2}+\frac 6{10}\,\mathbf{q}^{(4)}_{3,1}+\frac 1{10}\,\mathbf{q}^{(4)}_{3,0}\\
  \mathbf{q}^{(3)}_{0,3}&= \frac 1{10}\,\mathbf{q}^{(4)}_{3,3}+\frac 6{10}\,\mathbf{q}^{(4)}_{3,2}+\frac 3{10}\mathbf{q}^{(4)}_{3,1}\\
  \mathbf{q}^{(3)}_{0,4}&= \frac 2{5}\,\mathbf{q}^{(4)}_{3,3}+\frac 3{5}\,\mathbf{q}^{(4)}_{3,2}\\
  \mathbf{q}^{(3)}_{0,5}&= \mathbf{q}^{(4)}_{3,3}.
\end{split}
\end{equation}

A conclusion from the result above is that the control points $\mathbf{q}^{(3)}_{i,j}$ for $i=2,\ldots,5,$ $j=0,1$, and $i=0,1,\, j=2,\ldots,5$, are uniquely defined. For $i,j=0,1$, on the other hand, the control points $\mathbf{q}^{(3)}_{i,j}$ are defined in two different ways above. Concerning $\mathbf{q}^{(3)}_{0,0}$ there is obviously no problem. The fact that we have used the compatibility conditions (\ref{4.14}) implies equality in the two different ways defining $\mathbf{q}^{(3)}_{1,0}$ and $\mathbf{q}^{(3)}_{0,1}$. The last control point $\mathbf{q}^{(3)}_{1,1}$ must be chosen in such a way that the relations (\ref{4.11}) are fulfilled and consequently it is uniquely defined. Therefore, let us calculate $\mathbf{q}^{(3)}_{1,1}$ in terms of $\{\mathbf{q}^{(1)}_{i,j}\}_{i,j=2}^3$ in the two different ways that this control point has been defined in (\ref{4.16}) and (\ref{4.17}) respectively, in order to see what conditions are needed for fulfilling uniqueness. First, it follows
\begin{equation*}
  \begin{split}
    Q_{2,3}=&25(\mathbf{q}^{(3)}_{1,1}-\mathbf{q}^{(3)}_{1,0}-\mathbf{q}^{(3)}_{0,1}+\mathbf{q}^{(3)}_{0,0})
    =9\lambda_{1,4}\,(\mathbf{q}^{(2)}_{1,3}-\mathbf{q}^{(2)}_{1,2})+6\alpha_{2,3}\,(\mathbf{q}^{(2)}_{0,3}-\mathbf{q}^{(2)}_{0,2})\\&
    + 9\beta_{2,3}^{(1)}\,(\mathbf{q}^{(2)}_{1,3}-\mathbf{q}^{(2)}_{0,3})
    -15\lambda_{1,4}\,(\mathbf{q}^{(2)}_{0,3}-\mathbf{q}^{(2)}_{0,2})\\
    =&9\lambda_{1,4}\,(\mathbf{q}^{(2)}_{1,3}-\mathbf{q}^{(2)}_{0,3}+\mathbf{q}^{(2)}_{0,3}-\mathbf{q}^{(2)}_{0,2}-(\mathbf{q}^{(2)}_{1,2}-\mathbf{q}^{(2)}_{0,2}))\\&
    +(6\alpha_{2,3}-15\lambda_{1,4})\,(\mathbf{q}^{(2)}_{0,3}-\mathbf{q}^{(2)}_{0,2})
    +9\beta_{2,3}^{(1)}\,(\mathbf{q}^{(2)}_{1,3}-\mathbf{q}^{(2)}_{0,3}) \\
    =&(9\lambda_{1,4}+9\beta_{2,3}^{(1)})\,\lambda_{1,2}(\mathbf{q}^{(1)}_{3,3}-\mathbf{q}^{(1)}_{2,3})
    +(6\alpha_{2,3}-6\lambda_{1,4})\,(\mathbf{q}^{(1)}_{3,3}-\mathbf{q}^{(1)}_{3,2})\\&
    -9\lambda_{1,4}\lambda_{1,2}(\mathbf{q}^{(1)}_{3,2}-\mathbf{q}^{(1)}_{2,2}).
  \end{split}
\end{equation*}
Second, we consider the same expression as above given in (\ref{4.17}). We have
\begin{equation*}
  \begin{split}
    Q_{4,3}=&25(\mathbf{q}^{(3)}_{1,1}-\mathbf{q}^{(3)}_{0,1}-\mathbf{q}^{(3)}_{1,0}+\mathbf{q}^{(3)}_{0,0})
    =9\lambda_{1,2}\,(\mathbf{q}^{(4)}_{3,1}-\mathbf{q}^{(4)}_{2,1})+6\alpha_{4,3}\,(\mathbf{q}^{(4)}_{3,0}-\mathbf{q}^{(4)}_{2,0})\\&
    + 9\beta_{4,3}^{(1)}\,(\mathbf{q}^{(4)}_{3,1}-\mathbf{q}^{(4)}_{3,0})
    -15\lambda_{1,2}\,(\mathbf{q}^{(4)}_{3,0}-\mathbf{q}^{(4)}_{2,0})\\
    =&9\lambda_{1,2}\,(\mathbf{q}^{(4)}_{3,1}-\mathbf{q}^{(4)}_{3,0}
    +\mathbf{q}^{(4)}_{3,0}-\mathbf{q}^{(4)}_{2,0}-(\mathbf{q}^{(4)}_{2,1}-\mathbf{q}^{(4)}_{2,0}))\\&
    +(6\alpha_{4,3}-15\lambda_{1,2})\,(\mathbf{q}^{(4)}_{3,0}-\mathbf{q}^{(4)}_{2,0})+9\beta_{4,3}^{(1)}\,(\mathbf{q}^{(4)}_{3,1}-\mathbf{q}^{(4)}_{3,0}) \\
    =&(9\lambda_{1,2}+9\beta_{4,3}^{(1)})\,\lambda_{1,4}(\mathbf{q}^{(1)}_{3,3}-\mathbf{q}^{(1)}_{3,2})
    +(6\alpha_{4,3}-6\lambda_{1,2})\,(\mathbf{q}^{(1)}_{3,3}-\mathbf{q}^{(1)}_{2,3})\\&
    -9\lambda_{1,2}\lambda_{1,4}(\mathbf{q}^{(1)}_{2,3}-\mathbf{q}^{(1)}_{2,2}).
  \end{split}
\end{equation*}
Finally, using (\ref{1.9}) we get
\begin{equation*}
  \begin{split}
    Q_{2,3}-Q_{4,3}=&\big(9\lambda_{1,2}\beta_{2,3}^{(1)}-6(\alpha_{4,3}-\lambda_{1,2})\big)\,(\mathbf{q}^{(1)}_{3,3}-\mathbf{q}^{(1)}_{2,3})\\&+\big(6(\alpha_{2,3}-\lambda_{1,4})-9\lambda_{1,4}\beta_{4,3}^{(1)}\big)\,(\mathbf{q}^{(1)}_{3,3}-\mathbf{q}^{(1)}_{3,2})\\&
    +9\lambda_{1,4}\lambda_{1,2}\,(\mathbf{q}^{(1)}_{3,3}-\mathbf{q}^{(1)}_{2,3})
    -9\lambda_{1,2}\lambda_{1,4}\,(\mathbf{q}^{(1)}_{3,3}-\mathbf{q}^{(1)}_{3,2})\\&
    -9\lambda_{1,4}\lambda_{1,2}\,(\mathbf{q}^{(1)}_{3,2}-\mathbf{q}^{(1)}_{2,2})
    +9\lambda_{1,2}\lambda_{1,4}\,(\mathbf{q}^{(1)}_{2,3}-\mathbf{q}^{(1)}_{2,2})\\
    =&\big(9\lambda_{1,2}\beta_{2,3}^{(1)}-6(\alpha_{4,3}-\lambda_{1,2})\big)\,(\mathbf{q}^{(1)}_{3,3}-\mathbf{q}^{(1)}_{2,3})\\&+\big(6(\alpha_{2,3}-\lambda_{1,4})-9\lambda_{1,4}\beta_{4,3}^{(1)}\big)\,(\mathbf{q}^{(1)}_{3,3}-\mathbf{q}^{(1)}_{3,2})=0.
  \end{split}
\end{equation*}
The fact that the vectors $\mathbf{q}^{(1)}_{3,3}-\mathbf{q}^{(1)}_{2,3}$ and $\mathbf{q}^{(1)}_{3,3}-\mathbf{q}^{(1)}_{3,2}$ are linearly independent implies that the control point $\mathbf{q}^{(3)}_{1,1}$ is uniquely defined if and only if the coefficients $\alpha_{i,3}$, $\beta_{i,3}^{(1)}$ for $i=2,4$, satisfy
\begin{equation}\label{4.19}
  \begin{split}
    &3\lambda_{1,2}\beta_{2,3}^{(1)}=2\alpha_{4,3}-2\lambda_{1,2}\\
    &3\lambda_{1,4}\beta_{4,3}^{(1)}=2\alpha_{2,3}-2\lambda_{1,4}.
  \end{split}
\end{equation}
With the above criteria fulfilled the control points $\mathbf{q}^{(3)}_{i,j}$ for $i=0,1,\ldots,5$, $j=0,1$, and for $i=0,1$, $j=0,1,\ldots,5$, are uniquely defined by (\ref{4.15})--(\ref{4.18}). Concerning the other control points $\mathbf{q}^{(3)}_{i,j}$,\, $i,j=2,3,4,5$, there is no criteria given here how to choose these. Use any suitable method to determine those control points. The same is true for the parameters $\beta^{(2)}_{i,3}$ with $i=2,4$, while there may be certain restrictions concerning $\lambda_{i,3}(1)$, $\kappa_{i,3}(1)$ for $i=2,4$.

In the case, referred to on page~13, where $\kappa_{1,2}(0)\ne0$ or $\kappa_{1,4}(0)\ne0$, then the coefficient relations in (\ref{4.19}) are replaced by
\begin{equation*}
  \begin{split}
    &3\lambda_{1,2}\beta_{2,3}^{(1)}=2\alpha_{4,3}-2\lambda_{1,2}-\lambda_{1,2}\,\kappa_{1,4}(0)\\
    &3\lambda_{1,4}\beta_{4,3}^{(1)}=2\alpha_{2,3}-2\lambda_{1,4}-\lambda_{1,4}\,\kappa_{1,2}(0).
  \end{split}
\end{equation*}

Our goal so far has been to create a 4$^{th}$ patch in a non-complete 4-patch surface in such a way that the new surface is $G^1$-regular in spite of the fact that the functions $\lambda_{2,3}$ and $\lambda_{4,3}$ are non-constant. %and the connection between the newly created patch and its surroundings has non-constant $\lambda$- and $\kappa$-functions, in order to make the connection more flexible. More precisely all the functions $\lambda_{i,3}$ and $\kappa_{i,3}$ in (\ref{4.11}) are non-constant.
We have a further goal with this construction, which will be seen in the next subsection, but for the just mentioned purpose we can in fact choose $\mathbf{r}^{(3)}$ of bi-degree (4,4). Let
%In principle we can use polynomial functions $\lambda_{i,3}$, $\kappa_{i,3}$ for $i=2,4$ of lower degree than those in (\ref{4.141}). Let
\begin{equation*}
  \begin{split}
    &\beta_{i,3}^{(1)}=\frac13(2\beta_{i,3}+\kappa_{i,3}(0))\\
    &\beta_{i,3}^{(2)}=\frac13(2\beta_{i,3}+\kappa_{i,3}(1))\\
    &\alpha_{i,3}=\frac12(\lambda_{i,3}(0)+\lambda_{i,3}(1)),
  \end{split}
\end{equation*}
for $i=2,4$. Then we get the polynomials
%which implies that the patch $\mathbf{r}^{(3)}$ can be of bi-degree (4,4) and
\begin{eqnarray*}
  \lambda_{2,3}(u)&=&\lambda_{1,4}\,(1-u)+\lambda_{2,3}(1)\,u\\
  \lambda_{4,3}(v)&=&\lambda_{1,2}\,(1-v)+\lambda_{4,3}(1)\,v\\
  \kappa_{2,3}(u)&=&2\,\beta_{2,3}\,(1-u)u+\kappa_{2,3}(1)\,u^2\\
%\noalign{and}
\kappa_{4,3}(v)&=&2\,\beta_{4,3}^{(1)}\,(1-v)v+\kappa_{4,3}(1)\,v^2.
\end{eqnarray*}
 In this case necessary and sufficient conditions in (\ref{4.19}) are replaced by the next conditions
\begin{equation}
  \begin{split}
    &2\lambda_{1,2}\beta_{2,3}=\lambda_{4,3}(1)-\lambda_{1,2}=\lambda_{4,3}(1)-\lambda_{4,3}(0)\\
    &2\lambda_{1,4}\beta_{4,3}=\lambda_{2,3}(1)-\lambda_{1,4}=\lambda_{2,3}(1)-\lambda_{2,3}(0),
  \end{split}\label{4.193}
\end{equation}
which combined with (\ref{4.14}) imply uniqueness of the control point $\mathbf{q}^{(3)}_{1,1}$. On the other hand, this later solution is not flexible enough in the sense that it can not be used in the next application. This follows from (\ref{4.193}) in the context of the following subsection.

In the next subsection we will use the result we just received in order to fill a hole in a surface.

\subsection{Filling a hole in a surface}
In this subsection we consider a surface with an interior hole, more specific an incomplete 9-patch surface. By that we mean a surface like the one in Figure~5 where the part denoted by $\mathbf{r}^{(5)}$ is not included in the surface. We assume that the incomplete 9-patch surface is $G^1$-regular. We also assume that each of the original 8 patches are represented by a Bezier polynomial of bi-degree $(3,3)$ and that along the boundary curve between any two of the patches holds that the function $\lambda$ is constant and the function $\kappa=0$. See (\ref{eqn4.01}). The problem to consider here is to create an interior patch $\mathbf{r}^{(5)}$ with as low bi-degree as possible keeping the $G^1$-regularity of the complete surface. This is the situation Sarraga considered in the paper \cite{sarraga}. With his assumption he needed bi-degree (6,6) for $\mathbf{r}^{(5)}$. In this subsection we construct a patch $\mathbf{r}^{(5)}$ of bi-degree (5,5).

In order to create a Bezier patch $\mathbf{r}^{(5)}$ of bi-degree $(5,5)$ we start by considering it as part of the 4-patch surface consisting of $\mathbf{r}^{(1)}$, $\mathbf{r}^{(4)}$, $\mathbf{r}^{(5)}$ and $\mathbf{r}^{(2)}$. From the assumption above we have, as in the previous subsection, that the functions $\lambda_{1,2}$, $\lambda_{1,4}$ are constant and   $\kappa_{1,2}=\kappa_{1,4}=0$. We know that there exist compatibility conditions which are necessary and sufficient in order to get a unique solution for the control points $\mathbf{q}^{(5)}_{i,j}$ for $i=0,1,\ldots,5,$ $j=0,1$, and $i=0,1,\, j=0,1,\ldots,5$. Next we do the same with $\mathbf{r}^{(5)}$, $\mathbf{r}^{(8)}$, $\mathbf{r}^{(9)}$ and $\mathbf{r}^{(6)}$.  We start by defining the missing patch $\mathbf{r}^{(5)}$ as part of the first 4-patch surface. We have to solve
\begin{equation}\label{4.201}
    \begin{split}
        \mathbf{r}^{(5)}_v(u,0) & = \lambda_{4,5}(u)\mathbf{r}^{(4)}_v(u,1)+ \kappa_{4,5}(u)\mathbf{r}^{(4)}_u(u,1)\\
        \mathbf{r}^{(5)}_u (0,v) & = \lambda_{2,5}(v)\mathbf{r}^{(2)}_u(1,v)+\kappa_{2,5}(v)\mathbf{r}^{(2)}_v(1,v),
    \end{split}
\end{equation}
and
\begin{equation}\label{4.202}
    \left . \begin{array}{l}
        \mathbf{r}^{(5)}(u,0) = \mathbf{r}^{(2)}(u,1)\\
        \mathbf{r}^{(5)}(0,v) = \mathbf{r}^{(4)}(1,v),
    \end{array}\right .
\end{equation}
with
\begin{equation}\label{4.203}
    \begin{array}{ll}
        &\lambda_{4,5}(u)=\lambda_{4,5}(0)\,(1-u)^2+2\,\alpha_{4,5}\,(1-u)u+\lambda_{4,5}(1)\,u^2\\
        &\lambda_{2,5}(v)=\lambda_{2,5}(0)\,(1-v)^2+2\,\alpha_{2,5}\,(1-v)v+\lambda_{2,5}(1)\,v^2\\
        &\kappa_{4,5}(u)=\kappa_{4,5}(0)\,(1-u)^3+3\,\beta_{4,5}^{(1)}\,(1-u)^2u+3\,\beta_{4,5}^{(2)}\,(1-u)u^2+\kappa_{4,5}(1)\,u^3\\
        &\kappa_{2,5}(v)=\kappa_{2,5}(0)\,(1-v)^3+3\,\beta_{2,5}^{(1)}\,(1-v)^2v+3\,\beta_{2,5}^{(2)}\,(1-v)v^2+\kappa_{2,5}(1)\,v^3,
    \end{array}
\end{equation}
where the parameters $u,v\in[0,1]$.

From the previous subsection we know that in order to get a unique solution for the control points $\mathbf{q}^{(5)}_{i,j}$ for $i=0,1,\ldots,5,$ $j=0,1$, and $i=0,1,\, j=0,1,\ldots,5$, necessary and sufficient conditions are
\begin{equation}
    \begin{array}{ll}
        \lambda_{4,5}(0)=\lambda_{1,2},\,\,&
        \lambda_{2,5}(0)=\lambda_{1,4}\\
        \kappa_{4,5}(0)=0,\,\,&
        \kappa_{2,5}(0)=0
    \end{array} \label{4.21}
\end{equation}
  and
\begin{equation}
    \begin{array}{l}
        3\,\lambda_{1,4}\,\beta_{4,5}^{(1)}=2\,(\alpha_{2,5}-\lambda_{1,4})\\
        3\,\lambda_{1,2}\,\beta_{2,5}^{(1)}=2\,(\alpha_{4,5}-\lambda_{1,2}).
    \end{array} \label{4.25}
\end{equation}\vspace{1cm}

\begin{figure}[hbtp]
  \begin{center}
    \leavevmode
    \includegraphics[width=7cm]{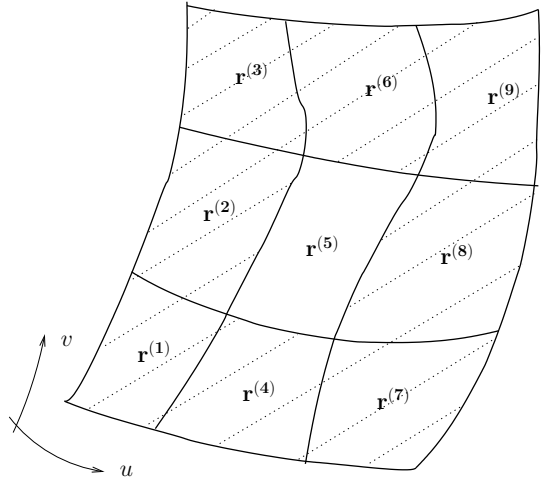}
    \caption{Filling a hole in a surface}
  \end{center}\label{fig4}
\end{figure}
As in the previous subsection, compare (\ref{4.15}) and (\ref{4.18}), we next get the solution of the boundary control points for the patch $\mathbf{r}^{(5)}$
\begin{equation}
    \begin{array}{l@{\hspace{1cm}}l}
        \mathbf{q}^{(5)}_{0,0}= \mathbf{q}^{(4)}_{0,3}&
        \mathbf{q}^{(5)}_{0,0}= \mathbf{q}^{(2)}_{3,0}\\
        \mathbf{q}^{(5)}_{1,0}= \frac35\,\mathbf{q}^{(4)}_{1,3}+\frac 25\,\mathbf{q}^{(4)}_{0,3}&
        \mathbf{q}^{(5)}_{0,1}= \frac35\,\mathbf{q}^{(2)}_{3,1}+\frac 25\,\mathbf{q}^{(2)}_{3,0}\\
  \mathbf{q}^{(5)}_{2,0}= \frac 3{10}\,\mathbf{q}^{(4)}_{2,3}+\frac 6{10}\,\mathbf{q}^{(4)}_{1,3}
  +\frac 1{10}\,\mathbf{q}^{(4)}_{0,3}&
  \mathbf{q}^{(5)}_{0,2}= \frac 3{10}\,\mathbf{q}^{(2)}_{3,2}+\frac 6{10}\,\mathbf{q}^{(2)}_{3,1}
  +\frac 1{10}\,\mathbf{q}^{(2)}_{3,0}\\
  \mathbf{q}^{(5)}_{3,0}= \frac 1{10}\,\mathbf{q}^{(4)}_{3,3}+\frac 6{10}\,\mathbf{q}^{(4)}_{2,3}
  +\frac 3{10}\,\mathbf{q}^{(4)}_{1,3}&
  \mathbf{q}^{(5)}_{0,3}= \frac 1{10}\,\mathbf{q}^{(2)}_{3,3}+\frac 6{10}\,\mathbf{q}^{(2)}_{3,2}
  +\frac 3{10}\,\mathbf{q}^{(2)}_{3,1}\\
  \mathbf{q}^{(5)}_{4,0}= \frac 25\,\mathbf{q}^{(4)}_{3,3}+\frac 35\,\mathbf{q}^{(4)}_{2,3}&
  \mathbf{q}^{(5)}_{0,4}= \frac 25\,\mathbf{q}^{(2)}_{3,3}+\frac 35\,\mathbf{q}^{(2)}_{3,2}\\
  \mathbf{q}^{(5)}_{5,0}= \mathbf{q}^{(4)}_{3,3}&
  \mathbf{q}^{(5)}_{0,5}= \mathbf{q}^{(2)}_{3,3}.
    \end{array}\label{4.22}
\end{equation}

We get, compare (\ref{4.16}) and (\ref{4.17}), part of its interior control points
\begin{equation}
  \begin{split}
    \mathbf{q}^{(5)}_{0,1}-\mathbf{q}^{(5)}_{0,0}=& \frac{3\,\lambda_{1,2}}5\,(\mathbf{q}^{(4)}_{0,3}-\mathbf{q}^{(4)}_{0,2})
    \\
    \mathbf{q}^{(5)}_{1,1}-\mathbf{q}^{(5)}_{1,0}=&\frac35\,\big(\frac{3\,\lambda_{1,2}}{5}\,(\mathbf{q}^{(4)}_{1,3}-\mathbf{q}^{(4)}_{1,2})+\frac{2\alpha_{4,5}}{5}\,(\mathbf{q}^{(4)}_{0,3}-\mathbf{q}^{(4)}_{0,2})
    +\frac{3\,\beta_{4,5}^{(1)}}{5}\,(\mathbf{q}^{(4)}_{1,3}-\mathbf{q}^{(4)}_{0,3})\big)
    \\
    \mathbf{q}^{(5)}_{2,1}-\mathbf{q}^{(5)}_{2,0}=&\frac35\,\big(\frac{3\,\lambda_{1,2}}{10}\,(\mathbf{q}^{(4)}_{2,3}-\mathbf{q}^{(4)}_{2,2})+\frac{6\,\alpha_{4,5}}{10}\,(\mathbf{q}^{(4)}_{1,3}-\mathbf{q}^{(4)}_{1,2})+\frac{\lambda_{7,8}}{10}\,(\mathbf{q}^{(4)}_{0,3}-\mathbf{q}^{(4)}_{0,2})\\&
    +\frac{6\,\beta_{4,5}^{(1)}}{10}\,(\mathbf{q}^{(4)}_{2,3}-\mathbf{q}^{(4)}_{1,3})+
    \frac{3\,\beta_{4,5}^{(2)}}{10}\,(\mathbf{q}^{(4)}_{1,3}-\mathbf{q}^{(4)}_{0,3})\big)
    \\
    \mathbf{q}^{(5)}_{3,1}-\mathbf{q}^{(5)}_{3,0}=& \frac35\,\big(\frac{\lambda_{1,2}}{10}\,(\mathbf{q}^{(4)}_{3,3}-\mathbf{q}^{(4)}_{3,2})+\frac{6\,\alpha_{4,5}}{10}\,(\mathbf{q}^{(4)}_{2,3}-\mathbf{q}^{(4)}_{2,2})+\frac{3\,\lambda_{7,8}}{10}\,(\mathbf{q}^{(4)}_{1,3}-\mathbf{q}^{(4)}_{1,2})\\&
    +\frac{3\,\beta_{4,5}^{(1)}}{10}\,(\mathbf{q}^{(4)}_{3,3}-\mathbf{q}^{(4)}_{2,3})+
    \frac{6\,\beta_{4,5}^{(2)}}{10}\,(\mathbf{q}^{(4)}_{2,3}-\mathbf{q}^{(4)}_{1,3})\big)
    \\
    \mathbf{q}^{(5)}_{4,1}-\mathbf{q}^{(5)}_{4,0}=&\frac35\,\big(\frac{2\,\alpha_{4,5}}{5}\,(\mathbf{q}^{(4)}_{3,3}-\mathbf{q}^{(4)}_{3,2})+\frac{3\,\lambda_{7,8}}{5}\,(\mathbf{q}^{(4)}_{2,3}-\mathbf{q}^{(4)}_{2,2})
+\frac{3\,\beta_{4,5}^{(2)}}{5}\,(\mathbf{q}^{(4)}_{3,3}-\mathbf{q}^{(4)}_{2,3})\big)
    \\
    \mathbf{q}^{(5)}_{5,1}-\mathbf{q}^{(5)}_{5,0}=&\frac{3\,\lambda_{7,8}}5\,(\mathbf{q}^{(4)}_{3,3}-\mathbf{q}^{(4)}_{3,2}),
\end{split}\label{4.23}
\end{equation}
and
\begin{equation}
  \begin{split}
    \mathbf{q}^{(5)}_{1,0}-\mathbf{q}^{(5)}_{0,0}=& \frac{3\,\lambda_{1,4}}5\,(\mathbf{q}^{(2)}_{3,0}-\mathbf{q}^{(2)}_{2,0})
    \\
    \mathbf{q}^{(5)}_{1,1}-\mathbf{q}^{(5)}_{0,1}=&\frac35\,\big(\frac{3\,\lambda_{1,4}}{5}\,(\mathbf{q}^{(2)}_{3,1}-\mathbf{q}^{(2)}_{2,1})+\frac{2\alpha_{2,5}}{5}\,(\mathbf{q}^{(2)}_{3,0}-\mathbf{q}^{(2)}_{2,0})
    +\frac{3\,\beta_{2,5}^{(1)}}{5}\,(\mathbf{q}^{(2)}_{3,1}-\mathbf{q}^{(2)}_{3,0})\big)
    \\
    \mathbf{q}^{(5)}_{1,2}-\mathbf{q}^{(5)}_{0,2}=&\frac35\,\big(\frac{3\,\lambda_{1,4}}{10}\,(\mathbf{q}^{(2)}_{3,2}-\mathbf{q}^{(2)}_{2,2})+\frac{6\,\alpha_{2,5}}{10}\,(\mathbf{q}^{(2)}_{3,1}-\mathbf{q}^{(2)}_{2,1})+\frac{\lambda_{3,6}}{10}\,(\mathbf{q}^{(2)}_{3,0}-\mathbf{q}^{(2)}_{2,0})\\&
    +\frac{6\,\beta_{2,5}^{(1)}}{10}\,(\mathbf{q}^{(2)}_{3,2}-\mathbf{q}^{(2)}_{3,1})+
    \frac{3\,\beta_{2,5}^{(2)}}{10}\,(\mathbf{q}^{(2)}_{3,1}-\mathbf{q}^{(2)}_{3,0})\big)
    \\
    \mathbf{q}^{(5)}_{1,3}-\mathbf{q}^{(5)}_{0,3}=& \frac35\,\big(\frac{\lambda_{1,4}}{10}\,(\mathbf{q}^{(2)}_{3,3}-\mathbf{q}^{(2)}_{2,3})+\frac{6\,\alpha_{2,5}}{10}\,(\mathbf{q}^{(2)}_{3,2}-\mathbf{q}^{(2)}_{2,2})+\frac{3\,\lambda_{3,6}}{10}\,(\mathbf{q}^{(2)}_{3,1}-\mathbf{q}^{(2)}_{2,1})\\&
    +\frac{3\,\beta_{2,5}^{(1)}}{10}\,(\mathbf{q}^{(2)}_{3,3}-\mathbf{q}^{(2)}_{3,2})+
    \frac{6\,\beta_{2,5}^{(2)}}{10}\,(\mathbf{q}^{(2)}_{3,2}-\mathbf{q}^{(2)}_{3,1})\big)
    \\
    \mathbf{q}^{(5)}_{1,4}-\mathbf{q}^{(5)}_{0,4}=&\frac35\,\big(\frac{2\,\alpha_{2,5}}{5}\,(\mathbf{q}^{(2)}_{3,3}-\mathbf{q}^{(2)}_{2,3})+\frac{3\,\lambda_{3,6}}{5}\,(\mathbf{q}^{(2)}_{3,2}-\mathbf{q}^{(2)}_{2,2})
+\frac{3\,\beta_{2,5}^{(2)}}{5}\,(\mathbf{q}^{(2)}_{3,3}-\mathbf{q}^{(2)}_{3,2})\big)
    \\
    \mathbf{q}^{(5)}_{1,5}-\mathbf{q}^{(5)}_{0,5}=&\frac{3\,\lambda_{3,6}}5\,(\mathbf{q}^{(2)}_{3,3}-\mathbf{q}^{(2)}_{2,3}),
\end{split}\label{4.24}
\end{equation}
where we also have used the fact that $\lambda_{4,5}(1)=\lambda_{7,8}$ and $\kappa_{4,5}(1)=0$ as well as $\lambda_{2,5}(1)=\lambda_{3,6}$ and $\kappa_{2,5}(1)=0$, which all follow from (\ref{4.2915}).

Let us continue with creating the next 4-patch surface containing the 4 patches $\mathbf{r}^{(5)}$, $\mathbf{r}^{(8)}$, $\mathbf{r}^{(9)}$ and $\mathbf{r}^{(6)}$, see Figure~5. With the natural change of indices we have here the same equations as in (\ref{4.201})--(\ref{4.202}) together with (\ref{4.203}). Thus, using the following compatibility conditions
\begin{equation}
    \begin{array}{ll}
        \lambda_{6,5}(1)=\lambda_{9,8},
        &\lambda_{8,5}(1)=\lambda_{9,6}\\
        \kappa_{6,5}(1)=0,
        &\kappa_{8,5}(1)=0,
    \end{array} \label{4.26}
\end{equation}
and
\begin{equation}
    \begin{array}{l}
        3\,\lambda_{9,6}\,\beta_{6,5}^{(2)}=2\,(\alpha_{8,5}-\lambda_{9,6})\\
        3\,\lambda_{9,8}\,\beta_{8,5}^{(2)}=2\,(\alpha_{6,5}-\lambda_{9,8}),
    \end{array} \label{4.291}
\end{equation}
combined with $\lambda_{6,5}(0)=\lambda_{3,2}$, $\kappa_{6,5}(0)=0$ and $\lambda_{8,5}(0)=\lambda_{7,4}$, $\kappa_{8,5}(0)=0$ from (\ref{4.2915}), we easily achieve the next control points, $\mathbf{q}^{(5)}_{i,j}$ for $i=0,1,\ldots,5,$ $j=4,5$, and $i=4,5,\, j=0,1,\ldots,5$, in a similar way as before. We have
\begin{equation}\label{4.27}
    \begin{array}{l@{\hspace{1cm}}l}
    \mathbf{q}^{(5)}_{0,5}= \mathbf{q}^{(6)}_{0,0}&
    \mathbf{q}^{(5)}_{5,0}= \mathbf{q}^{(8)}_{0,0}\\
    \mathbf{q}^{(5)}_{1,5}= \frac35\,\mathbf{q}^{(6)}_{1,0}+\frac 25\,\mathbf{q}^{(6)}_{0,0}&
    \mathbf{q}^{(5)}_{5,1}= \frac35\,\mathbf{q}^{(8)}_{0,1}+\frac 25\,\mathbf{q}^{(8)}_{0,0}\\
    \mathbf{q}^{(5)}_{2,5}= \frac 3{10}\,\mathbf{q}^{(6)}_{2,0}+\frac 6{10}\,\mathbf{q}^{(6)}_{1,0}
    +\frac 1{10}\,\mathbf{q}^{(6)}_{0,0}&
    \mathbf{q}^{(5)}_{5,2}= \frac 3{10}\,\mathbf{q}^{(8)}_{0,2}+\frac 6{10}\,\mathbf{q}^{(8)}_{0,1}
    +\frac 1{10}\,\mathbf{q}^{(8)}_{0,0}\\
    \mathbf{q}^{(5)}_{3,5}= \frac 1{10}\,\mathbf{q}^{(6)}_{3,0}+\frac 6{10}\,\mathbf{q}^{(6)}_{2,0}
    +\frac 3{10}\,\mathbf{q}^{(6)}_{1,0}&
    \mathbf{q}^{(5)}_{5,3}= \frac 1{10}\,\mathbf{q}^{(8)}_{0,3}+\frac 6{10}\,\mathbf{q}^{(8)}_{0,2}
    +\frac 3{10}\,\mathbf{q}^{(8)}_{0,1}\\
    \mathbf{q}^{(5)}_{4,5}= \frac 25\mathbf{q}^{(6)}_{3,0}+\frac 35\mathbf{q}^{(6)}_{2,0}&
    \mathbf{q}^{(5)}_{5,4}= \frac 25\,\mathbf{q}^{(8)}_{0,3}+\frac 35\,\mathbf{q}^{(8)}_{0,2}\\
    \mathbf{q}^{(5)}_{5,5}= \mathbf{q}^{(6)}_{3,0}&
    \mathbf{q}^{(5)}_{5,5}= \mathbf{q}^{(8)}_{0,3},
    \end{array}
\end{equation}
and
\begin{equation}\label{4.28}
  \begin{split}
    \mathbf{q}^{(5)}_{0,4}-\mathbf{q}^{(5)}_{0,5}=& \frac{3\,\lambda_{3,2}}5\,(\mathbf{q}^{(6)}_{0,0}-\mathbf{q}^{(6)}_{0,1})
    \\
    \mathbf{q}^{(5)}_{1,4}-\mathbf{q}^{(5)}_{1,5}=&
    \frac35\,\big(\frac{3\,\lambda_{3,2}}{5}\,(\mathbf{q}^{(6)}_{1,0}-\mathbf{q}^{(6)}_{1,1})
    +\frac{2\,\alpha_{6,5}}{5}\,(\mathbf{q}^{(6)}_{0,0}-\mathbf{q}^{(6)}_{0,1})
    +\frac{3\,\beta^{(1)}_{6,5}}{5}\,(\mathbf{q}^{(6)}_{1,0}-\mathbf{q}^{(6)}_{0,0})\big)
    \\
    \mathbf{q}^{(5)}_{2,4}-\mathbf{q}^{(5)}_{2,5}=&
    \frac35\,\big(\frac{3\,\lambda_{3,2}}{10}\,(\mathbf{q}^{(6)}_{2,0}-\mathbf{q}^{(6)}_{2,1})
    +\frac{6\,\alpha_{6,5}}{10}\,(\mathbf{q}^{(6)}_{1,0}-\mathbf{q}^{(6)}_{1,1})
    +\frac{\lambda_{9,8}}{10}\,(\mathbf{q}^{(6)}_{0,0}-\mathbf{q}^{(6)}_{0,1})\\&\makebox[2mm]{}
    +\frac{6\,\beta^{(1)}_{6,5}}{10}\,(\mathbf{q}^{(6)}_{2,0}-\mathbf{q}^{(6)}_{1,0})
    +\frac{3\,\beta^{(2)}_{6,5}}{10}\,(\mathbf{q}^{(6)}_{1,0}-\mathbf{q}^{(6)}_{0,0})\big)
    \\
    \mathbf{q}^{(5)}_{3,4}-\mathbf{q}^{(5)}_{3,5}=& \frac35\,\big(\frac{\lambda_{3,2}}{10}\,(\mathbf{q}^{(6)}_{3,0}-\mathbf{q}^{(6)}_{3,1})
    +\frac{6\,\alpha_{6,5}}{10}\,(\mathbf{q}^{(6)}_{2,0}-\mathbf{q}^{(6)}_{2,1})
    +\frac{3\,\lambda_{9,8}}{10}\,(\mathbf{q}^{(6)}_{1,0}-\mathbf{q}^{(6)}_{1,1})\\&\makebox[2mm]{}
    +\frac{3\,\beta^{(1)}_{6,5}}{10}\,(\mathbf{q}^{(6)}_{3,0}-\mathbf{q}^{(6)}_{2,0})+
    \frac{6\,\beta^{(2)}_{6,5}}{10}\,(\mathbf{q}^{(6)}_{2,0}-\mathbf{q}^{(6)}_{1,0})\big)
    \\
    \mathbf{q}^{(5)}_{4,4}-\mathbf{q}^{(5)}_{4,5}=&
    \frac35\,\big(\frac{2\,\alpha_{6,5}}5\,(\mathbf{q}^{(6)}_{3,0}-\mathbf{q}^{(6)}_{3,1})
    +\frac{3\,\lambda_{9,8}}5\,(\mathbf{q}^{(6)}_{2,0}-\mathbf{q}^{(6)}_{2,1})
    +\frac{3\,\beta^{(2)}_{6,5}}5\,(\mathbf{q}^{(6)}_{3,0}-\mathbf{q}^{(6)}_{2,0})\big)
    \\
    \mathbf{q}^{(5)}_{5,4}-\mathbf{q}^{(5)}_{5,5}=&
    \frac{3\,\lambda_{9,8}}5\,(\mathbf{q}^{(6)}_{3,0}-\mathbf{q}^{(6)}_{3,1})
  \end{split}
\end{equation}
and
\begin{equation}\label{4.29}
  \begin{split}
    \mathbf{q}^{(5)}_{4,0}-\mathbf{q}^{(5)}_{5,0}=& \frac{3\,\lambda_{7,4}}5\,(\mathbf{q}^{(8)}_{0,0}-\mathbf{q}^{(8)}_{1,0})
    \\
    \mathbf{q}^{(5)}_{4,1}-\mathbf{q}^{(5)}_{5,1}=&
    \frac35\,\big(\frac{3\,\lambda_{7,4}}{5}\,(\mathbf{q}^{(8)}_{0,1}-\mathbf{q}^{(8)}_{1,1})
    +\frac{2\,\alpha_{8,5}}{5}\,(\mathbf{q}^{(8)}_{0,0}-\mathbf{q}^{(8)}_{1,0})
    +\frac{3\,\beta^{(1)}_{8,5}}{5}\,(\mathbf{q}^{(8)}_{0,1}-\mathbf{q}^{(8)}_{0,0})\big)
    \\
    \mathbf{q}^{(5)}_{4,2}-\mathbf{q}^{(5)}_{5,2}=&
    \frac35\,\big(\frac{3\,\lambda_{7,4}}{10}\,(\mathbf{q}^{(8)}_{0,2}-\mathbf{q}^{(8)}_{1,2})
    +\frac{6\,\alpha_{8,5}}{10}\,(\mathbf{q}^{(8)}_{0,1}-\mathbf{q}^{(8)}_{1,1})
    +\frac{\lambda_{9,6}}{10}\,(\mathbf{q}^{(8)}_{0,0}-\mathbf{q}^{(8)}_{1,0})\\&\makebox[2mm]{}
    +\frac{6\,\beta^{(1)}_{8,5}}{10}\,(\mathbf{q}^{(8)}_{0,2}-\mathbf{q}^{(8)}_{0,1})
    +\frac{3\,\beta^{(2)}_{8,5}}{10}\,(\mathbf{q}^{(8)}_{0,1}-\mathbf{q}^{(8)}_{0,0})\big)
    \\
    \mathbf{q}^{(5)}_{4,3}-\mathbf{q}^{(5)}_{5,3}=& \frac35\,\big(\frac{\lambda_{7,4}}{10}\,(\mathbf{q}^{(8)}_{0,3}-\mathbf{q}^{(8)}_{1,3})
    +\frac{6\,\alpha_{8,5}}{10}\,(\mathbf{q}^{(8)}_{0,2}-\mathbf{q}^{(8)}_{1,2})
    +\frac{3\,\lambda_{9,6}}{10}\,(\mathbf{q}^{(8)}_{0,1}-\mathbf{q}^{(8)}_{1,1})\\&\makebox[2mm]{}
    +\frac{3\,\beta^{(1)}_{8,5}}{10}\,(\mathbf{q}^{(8)}_{0,3}-\mathbf{q}^{(8)}_{0,2})+
    \frac{6\,\beta^{(2)}_{8,5}}{10}\,(\mathbf{q}^{(8)}_{0,2}-\mathbf{q}^{(8)}_{0,1})\big)
    \\
    \mathbf{q}^{(5)}_{4,4}-\mathbf{q}^{(5)}_{5,4}=& \frac35\,\big(\frac{2\,\alpha_{8,5}}{5}\,(\mathbf{q}^{(8)}_{0,3}-\mathbf{q}^{(8)}_{1,3})
    +\frac{3\,\lambda_{9,6}}{5}\,(\mathbf{q}^{(8)}_{0,2}-\mathbf{q}^{(8)}_{1,2})
    +\frac{3\,\beta^{(2)}_{8,5}}{5}\,(\mathbf{q}^{(8)}_{0,3}-\mathbf{q}^{(8)}_{0,2})\big)
    \\
    \mathbf{q}^{(5)}_{4,5}-\mathbf{q}^{(5)}_{5,5}=&
    \frac{3\,\lambda_{9,6}}5\,(\mathbf{q}^{(8)}_{0,3}-\mathbf{q}^{(8)}_{1,3}).
\end{split}
\end{equation}

When collecting all the above control point we may have lost uniqueness of the doubly defined control points $\mathbf{q}^{(5)}_{ij}$ for $i=4,5,\,\,j=0,1$, and $i=0,1,\,\,j=4,5$. Considering a 4-patch surface around each of the other two vertices of the patch $\mathbf{r}^{(5)}$, i.e., the two 4-patch surfaces constituting of $\mathbf{r}^{(2)}$, $\mathbf{r}^{(5)}$, $\mathbf{r}^{(6)}$, $\mathbf{r}^{(3)}$ and $\mathbf{r}^{(4)}$, $\mathbf{r}^{(7)}$, $\mathbf{r}^{(8)}$, $\mathbf{r}^{(5)}$ respectively, we have necessary and sufficient conditions for uniqueness of the just mentioned control points. These conditions are
\begin{equation}
    \begin{array}{ll}
        \lambda_{4,5}(1)=\lambda_{7,8},\,\,&\lambda_{8,5}(0)=\lambda_{7,4}\\
        \kappa_{4,5}(1)=0,&\kappa_{8,5}(0)=0\\
        \lambda_{6,5}(0)=\lambda_{3,2},&
        \lambda_{2,5}(1)=\lambda_{3,6}\\
        \kappa_{6,5}(0)=0,&\kappa_{2,5}(1)=0
    \end{array} \label{4.2915}
\end{equation}
and
\begin{equation}
    \begin{array}{l}
        3\,\lambda_{7,4}\,\beta_{4,5}^{(2)}=2\,(\alpha_{8,5}-\lambda_{7,4})\\
        3\,\lambda_{7,8}\,\beta_{8,5}^{(1)}=2\,(\alpha_{4,5}-\lambda_{7,8}),
    \end{array} \label{4.292}
\end{equation}
\begin{equation}
    \begin{array}{l}
        3\,\lambda_{3,6}\,\beta_{6,5}^{(1)}=2\,(\alpha_{2,5}-\lambda_{3,6})\\
        3\,\lambda_{3,2}\,\beta_{2,5}^{(2)}=2\,(\alpha_{6,5}-\lambda_{3,2}).
    \end{array} \label{4.293}
\end{equation}

First we need to fulfil the compatibility conditions (\ref{4.25}), (\ref{4.291}), (\ref{4.292}) and (\ref{4.293}), where we have four degree of freedom for the parameters $\alpha_{i,5}$, $\beta_{i,5}^{(1)}$ and $\beta_{i,5}^{(2)}$ with $i=2,4,6,8$. After decided the value of the parameters, we have
partly defined an interior patch $\mathbf{r}^{(5)}$, where the control points are uniquely defined except for the undefined control points $\mathbf{q}^{(5)}_{ij}$ with $i,j=2,3$. Those points must be chosen in some way. One way to chose the undefined interior control points is by defining a Coons' patch from the boundary curves combined with their derivatives, see \cite{sarraga}, or use the next definition
\begin{equation}
    \begin{array}{l}
        \mathbf{q}^{(5)}_{2,2}=\mathbf{q}^{(5)}_{2,1}+\mathbf{q}^{(5)}_{1,2}-\mathbf{q}^{(5)}_{1,1}\\
        \mathbf{q}^{(5)}_{3,2}=\mathbf{q}^{(5)}_{3,1}+\mathbf{q}^{(5)}_{4,2}-\mathbf{q}^{(5)}_{4,1}\\
        \mathbf{q}^{(5)}_{2,3}=\mathbf{q}^{(5)}_{2,4}+\mathbf{q}^{(5)}_{1,3}-\mathbf{q}^{(5)}_{1,4}\\
        \mathbf{q}^{(5)}_{3,3}=\mathbf{q}^{(5)}_{3,4}+\mathbf{q}^{(5)}_{4,3}-\mathbf{q}^{(5)}_{4,4}.
    \end{array} \label{4.294}
\end{equation}

Thus, we have proved that it is possible to represent the interior patch $\mathbf{r}^{(5)}$ as a Bezier patch of bi-degree (5,5) in such a way that the complete 9-patch surface is $G^1$-continuous. On the other hand, the representation of $\mathbf{r}^{(5)}$ with non-constant $\lambda_{i,5}$-functions given above is very dependent of non-trivial $\kappa_{i,5}$-functions, because, if, on the contrary, the $\kappa_{i,5}$-functions are identically zero, then the $\lambda_{i,5}$-functions must be identically constant. This follows trivially from the compatibility conditions.

Suppose that the $\kappa_{i,5}$-functions are identically zero, but the $\lambda_{i,5}$-functions are non-constant, then the patch $\mathbf{r}^{(5)}$ must be of at least bi-degree (6,6). See also Sarraga \cite{sarraga}. In this case $\lambda_{4,5}$ and $\lambda_{2,5}$ are defined as
\begin{eqnarray*}
    &\lambda_{4,5}(u)=\lambda_{1,2}\,(1-u)^3+3\,\alpha_{4,5}^{(1)}\,(1-u)^2u+3\,\alpha_{4,5}^{(2)}\,(1-u)u^2
  +\lambda_{7,8}\,u^3\phantom{.}\\
    &\lambda_{2,5}(v)=\lambda_{1,4}\,(1-v)^3+3\,\alpha_{2,5}^{(1)}\,(1-v)^2v+3\,\alpha_{2,5}^{(2)}\,(1-u)u^2
  +\lambda_{3,6}\,v^3.
\end{eqnarray*}
The compatibility conditions will now be
\begin{eqnarray*}
    \alpha^{(1)}_{4,5}=\lambda_{1,2}\phantom{.}\\
    \alpha^{(1)}_{2,5}=\lambda_{1,4}.
\end{eqnarray*}
Continuing around the other vertices in the patch $\mathbf{r}^{(5)}$ we also get the following conditions
\begin{equation*}
    \begin{array}{lll}
        \alpha^{(2)}_{6,5}=\lambda_{9,8},&
        \alpha^{(2)}_{4,5}=\lambda_{7,8},&
        \alpha^{(1)}_{6,5}=\lambda_{3,2}\\
        \alpha^{(2)}_{8,5}=\lambda_{9,6},&
        \alpha^{(1)}_{8,5}=\lambda_{7,4},&
        \alpha^{(2)}_{2,5}=\lambda_{3,6}.
    \end{array}
\end{equation*}
The above conditions must be complemented with (\ref{4.21}), (\ref{4.26}) and (\ref{4.2915}). With these restrictions we get the unique control points of $\mathbf{r}^{(5)}$ except for the undefined points $\mathbf{q}^{(5)}_{i,j}$, $i,j=2,3,4$. Those interior points can be defined through a Coons' patch as above or in a similar way as in (\ref{4.294}), i.e.,
\begin{equation*}
    \begin{array}{l}
        \mathbf{q}^{(5)}_{2,2}=\mathbf{q}^{(5)}_{2,1}+\mathbf{q}^{(5)}_{1,2}-\mathbf{q}^{(5)}_{1,1}\\
        \mathbf{q}^{(5)}_{4,2}=\mathbf{q}^{(5)}_{4,1}+\mathbf{q}^{(5)}_{5,2}-\mathbf{q}^{(5)}_{5,1}\\
        \mathbf{q}^{(5)}_{2,4}=\mathbf{q}^{(5)}_{2,5}+\mathbf{q}^{(5)}_{1,4}-\mathbf{q}^{(5)}_{1,5}\\
        \mathbf{q}^{(5)}_{4,4}=\mathbf{q}^{(5)}_{4,5}+\mathbf{q}^{(5)}_{5,4}-\mathbf{q}^{(5)}_{5,5}
    \end{array} %\label{4.294}
\end{equation*}
and
\begin{equation*}
    \begin{array}{l}
        \mathbf{q}^{(5)}_{2,3}=\mathbf{q}^{(5)}_{1,3}+\frac12(\mathbf{q}^{(5)}_{2,2}+\mathbf{q}^{(5)}_{2,4})
            -\frac12(\mathbf{q}^{(5)}_{1,2}+\mathbf{q}^{(5)}_{1,4})
            =\mathbf{q}^{(5)}_{1,3}+\frac12(\mathbf{q}^{(5)}_{2,1}+\mathbf{q}^{(5)}_{2,5})
            -\frac12(\mathbf{q}^{(5)}_{1,1}+\mathbf{q}^{(5)}_{1,5})\\
        \mathbf{q}^{(5)}_{3,2}=\mathbf{q}^{(5)}_{3,1}+\frac12(\mathbf{q}^{(5)}_{2,2}+\mathbf{q}^{(5)}_{4,2})
            -\frac12(\mathbf{q}^{(5)}_{2,1}+\mathbf{q}^{(5)}_{4,1})
            =\mathbf{q}^{(5)}_{3,1}+\frac12(\mathbf{q}^{(5)}_{1,2}+\mathbf{q}^{(5)}_{5,2})
            -\frac12(\mathbf{q}^{(5)}_{1,1}+\mathbf{q}^{(5)}_{5,1})\\
        \mathbf{q}^{(5)}_{3,4}=\mathbf{q}^{(5)}_{3,5}+\frac12(\mathbf{q}^{(5)}_{2,4}+\mathbf{q}^{(5)}_{4,4})
            -\frac12(\mathbf{q}^{(5)}_{2,5}+\mathbf{q}^{(5)}_{4,5})
            =\mathbf{q}^{(5)}_{3,5}+\frac12(\mathbf{q}^{(5)}_{1,4}+\mathbf{q}^{(5)}_{5,4})
            -\frac12(\mathbf{q}^{(5)}_{1,5}+\mathbf{q}^{(5)}_{5,5})\\
        \mathbf{q}^{(5)}_{4,3}=\mathbf{q}^{(5)}_{5,3}+\frac12(\mathbf{q}^{(5)}_{4,2}+\mathbf{q}^{(5)}_{4,4})
            -\frac12(\mathbf{q}^{(5)}_{5,2}+\mathbf{q}^{(5)}_{5,4})
            =\mathbf{q}^{(5)}_{5,3}+\frac12(\mathbf{q}^{(5)}_{4,1}+\mathbf{q}^{(5)}_{4,5})
            -\frac12(\mathbf{q}^{(5)}_{5,1}+\mathbf{q}^{(5)}_{5,5})\\
            \end{array} %\label{4.294}
\end{equation*}
and finally
\begin{equation*}
    \begin{array}{l}
        \mathbf{q}^{(5)}_{3,3}=\frac12(\mathbf{q}^{(5)}_{2,3}+\mathbf{q}^{(5)}_{4,3}+\mathbf{q}^{(5)}_{3,2}
            +\mathbf{q}^{(5)}_{3,4})-\frac14(\mathbf{q}^{(5)}_{2,2}+\mathbf{q}^{(5)}_{4,2}
            +\mathbf{q}^{(5)}_{2,4}+\mathbf{q}^{(5)}_{4,2}).\\
            \end{array} %\label{4.294}
\end{equation*}

Observe, the parameters $\alpha^{(j)}_{i,5}$ for $i=2,4,6,8,$ and $j=1,2,$ are here completely decided from the compatibility conditions.

A further observation is that if the patch $\mathbf{r}^{(5)}$ is of bi-degree (4,4) the compatibility conditions, see (\ref{4.193}), in this case are not flexible enough to solve our problem. %In fact, either the polynomial functions $\kappa_{i,5}$$\lambda_{i,5}$ of degree one and  of degree degenerate to be identically constant and zero respectively.

We will use the main result of this subsection in the next part, where we will study a way to create fillet surfaces. \vspace{0.5cm}

\begin{figure}[h]\label{fig5}
  \begin{center}
    \leavevmode
     \includegraphics[width=8cm]{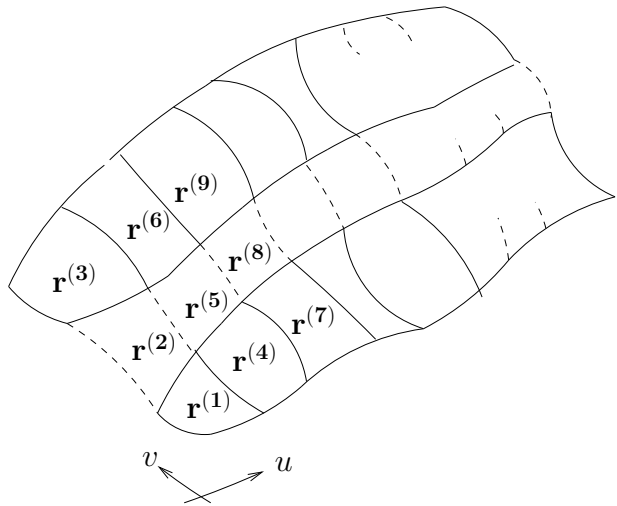}
    \caption{Creating a fillet}
  \end{center}
\end{figure}

\subsection{Creating a fillet surface}
 Consider two surfaces of regularity $G^1$, each consisting of the patches $\mathbf{r}^{(1+3n)}$ for $n=0,1,\ldots,N-1$ and $\mathbf{r}^{(3+3n)}$ for $n=0,1,\ldots,N-1$. We want to connect these two surfaces by creating a fillet surface in such a way that altogether there will be one complete $G^1$-surface. We assume that all the patches so far are of bi-degree (3,3) and that the connection between any two patches satisfies that the $\lambda$-function is constant and the $\kappa$-function is identically zero. See Figure~6.

Our first step is to create the patches $\mathbf{r}^{(2+6n)}$ for
$n=0,1,\ldots,[(N+1)/2]-1$, where $[N/2]$ denotes the integer part of
$N/2$. Let the patch $\mathbf{r}^{(2+6n)}$ be defined as a Bezier
patch of bi-degree $(3,3)$ connecting $\mathbf{r}^{(1+6n)}$ and
$\mathbf{r}^{(3+6n)}$ in a $G^1$-regular way. Furthermore, there are
no non-zero $\kappa$-function towards the two neighboring patches. In
the next step we create the patches in between, i.e.,
$\mathbf{r}^{(5+6n)}$ for $n=0,1,\ldots,[N/2]-1$, as we did in the
previous subsection. All together, this completes the construction of
the fillet surface.\vspace{5mm}

\noindent
{\bf Acknowledgement.} I like to thank Roger Andersson for valuable discussions.

%This patch we do construct as the missing interior patch in the uncomplete 9-patch surface as in the previous subsection.

%\newpage

  \endlines{Department of Mathematical Sciences\cr Chalmers University of Technology and University of Gothenburg\cr
    SE--412\thinspace96 Gothenburg, Sweden\cr
   bo@chalmers.se\cr
        }
\end{document}